\documentclass[preprint,12pt,english]{elsarticle}

\usepackage{fullpage}
\usepackage{amsmath,amsfonts,amsthm,amssymb}
\usepackage{color,xcolor}
\usepackage{ifpdf}
\usepackage{psfrag}
\usepackage{graphicx,graphics}
\usepackage{hhtensor}
\usepackage{comment}
\usepackage[small]{caption}
\usepackage{subfigure}
\usepackage{tikz}
\usepackage{mathtools}
\usepackage{bbm}
\usepackage{algorithm,algorithmic}
\usepackage{pgfplots}
\usepackage{pgfplotstable}
\pgfplotsset{compat = 1.3}
\usetikzlibrary{external}

\usepackage{xparse}
\usepackage{bigints}

\usepackage{url}
\usepackage{hyperref}

\usepackage[capitalise,noabbrev]{cleveref}
\usepackage{todonotes}
\newtheorem{theorem}{Theorem}[section]

\newtheorem{proposition}[theorem]{Proposition}
\newtheorem{remark}[theorem]{Remark}

\usepackage{cancel}
\usepackage{algorithm,algorithmic}
\usepackage{enumitem}
\newcommand{\diff}[1]{{\mathrm{d}{#1}}}
\newcommand{\bu}{\mathbf{U}}
\newcommand{\bF}{\mathbf{F}}
\newcommand{\bG}{\mathbf{G}}
\newcommand{\bU}{\mathbf{U}}
\newcommand{\bH}{\mathbf{H}}
\newcommand{\bA}{\mathbf{A}}
\newcommand{\bR}{\mathbf{R}}
\newcommand{\bB}{\mathbf{B}}

\newcommand{\xip}{x_{i+1/2}}
\newcommand{\xin}{x_{i-1/2}}
\newcommand{\iip}{{i+1/2}}
\newcommand{\iin}{{i-1/2}}
\newcommand{\bS}{\mathbf{S}}
\usepackage{bbm}

\def\R{\mathbb{R}}

\definecolor{darkspringgreen}{rgb}{0., 0.55, 0.3}
\definecolor{dartmouthgreen}{rgb}{0.05, 0.5, 0.06}
\definecolor{etonblue}{rgb}{0.59, 0.78, 0.64}
\definecolor{airforceblue}{rgb}{0., 0.4, 0.66}
\definecolor{arylideyellow}{rgb}{0.91, 0.84, 0.42}
\definecolor{emerald}{rgb}{0.31, 0.78, 0.47}
\definecolor{uclagold}{rgb}{1.0, 0.7, 0.0}
\definecolor{cadmiumorange}{rgb}{0.93, 0.53, 0.18}
\definecolor{redcomment}{rgb}{1.0, 0, 0.0}

\newsavebox{\DelimiterBox}
\newlength{\DelimiterHeight}
\newlength{\DelimiterDepth}
\newsavebox{\ArgumentBox}
\newlength{\ArgumentHeight}
\newlength{\ArgumentDepth}
\newlength{\ResizedDelimiterHeight}
\newlength{\ResizedDelimiterDepth}

\ifluatex

\else

\fi

\hypersetup{
pdftitle={}
pdfauthor={},
pdfpagemode=UseOutlines,
linkbordercolor=0 0 0,
linkcolor=red,
citecolor=blue,
colorlinks=true
}

\journal{Elsevier}

\begin{document}

\begin{frontmatter}


\title{High order global flux schemes for general steady state preservation of shallow water moment equations with non-conservative products}

\author[lab1]{Mirco Ciallella}
\author[lab2,lab3]{Julian Koellermeier}

\address[lab1]{Laboratoire Jacques-Louis Lions, Universit\'e Paris Cit\'e, CNRS, UMR 7598, Paris, France}
\address[lab2]{Department of Mathematics, Computer Science and Statistics, Ghent University, Ghent, Belgium}
\address[lab3]{Bernoulli Institute, University of Groningen, Groningen, the Netherlands}

\begin{abstract}
Shallow water moment equations are reduced-order models for free-surface flows that allow to represent vertical variations of the velocity profile at the expense of additional evolution equations for a number of additional variables, so called moments. This introduces non-linear non-conservative products in the system, which make the analytical characterization of steady states much harder if not impossible.
The lack of analytical steady states poses a challenge for the design of well-balanced schemes, which aim at preserving such steady states as crucial in many applications.

In this work, we present a family of fully well-balanced, high-order WENO finite volume methods for general hyperbolic balance laws with non-conservative products like the shallow water moment equations, for which no analytical steady states are available.
The schemes are based on the flux globalization approach, in which both source terms and non-conservative products are integrated with a tailored high order quadrature in the divergence term. The resulting
global flux is then reconstructed instead of the conservative variables to preserve all steady states.
Numerical tests show the optimal convergence of the method and a significant error reduction for steady state solutions. Furthermore, we provide a numerical comparison of perturbed steady states for different families of shallow water moment equations, which illustrates the flexibility of our method that is valid for general equations without prior knowledge of steady states.
\end{abstract}



\end{frontmatter}

\section{Introduction}

Free-surface flows are widely researched, with applications ranging from
tsunami modeling to river estuary and flood simulations \cite{bi2015mixed,ciallella2025high,macias2020performance}.
Although the underlying incompressible Navier-Stokes equations accurately model
these flows, the computational cost of such simulations is often either prohibitive, e.g.,
when dealing with large scale problems, or unnecessary, e.g., in simplified shallow water conditions.
The shallow water equations (SWE) provide a reduced model obtained by depth-averaging the
incompressible Navier-Stokes equations, assuming a constant velocity profile in the vertical direction.
However, this assumption may lead to important inaccuracies in the prediction of flows
with vertical variations in the velocity profile, e.g., due to friction terms or in the presence of physical phenomena such as sediment transport \cite{garres2020shallow}. But already in tsunami or dam break simulations,
this assumption is often violated \cite{koellermeier2020analysis}.

As an alternative to multilayer models \cite{kim2009two},
the shallow water moment equations (SWME) have been developed starting with \cite{kowalski2019moment}.
These models are based on an expansion of the velocity profile in vertical direction using
orthogonal Legendre polynomials. The system is closed by deriving additional evolution equations for the expansion coefficients, the so-called moments.
The resulting models are more accurate the more moments are included in the expansion,
but they also lack hyperbolicity, a necessary condition for numerical stability, as discussed in detail in \cite{koellermeier2020analysis}.
In the same work, a new family of hyperbolic models called hyperbolic shallow water equations (HSWME)
was derived based on a regularization of the original SWME model.
The numerical experiments in \cite{koellermeier2020analysis,kowalski2019moment} show good model accuracy for tests including friction terms, while only a flat bottom topography was considered.

As studied for the classical SWE \cite{alcrudo2001exact,bernetti2008exact}, it becomes fundamental
for practical applications to include a non-flat bottom topography in the SWME.
The addition of a varying bathymetry to the SWME leads to an intricate system which admits a large class of
steady states, given by the interaction between the flux and the forcing terms (including friction and bathymetry).
Since many practical applications can be modeled as perturbations of these steady states,
it becomes crucial to avoid numerical artifacts in simulations. This requires to design numerical schemes that preserve such steady states at the discrete level.

For simple models like the SWE, this typically follows a two-step procedure: 1) study the steady states of the models, and 2) design tailored well-balanced schemes \cite{bermudez1994upwind}, which preserve the analytically computed
steady states by balancing the flux and the source terms at the discrete level such that they cancel out.
The design of well-balanced schemes is a very active research topic in the literature
\cite{audusse2004fast,noelle2006well,gallardo2007well,canestrelli2009well,diaz2013high,bollermann2013well,
berthon2016fully,arpaia2018r,castro2020well,micalizzi2024novel}
and it is out of scope to give a complete overview here.

Recently, there has been increased interest in developing well-balanced method for the new
SWME \cite{koellermeier2022steady}. However, it was shown that the first step mentioned above,
i.e.\ the analytical characterization of the steady states, can be performed only for the first
order moment model, where the velocity profile is linear and only one additional moment needs to be taken into account. This limitation poses several
problems to the design and generalization of well-balanced schemes. To overcome this issue
the authors in \cite{koellermeier2022steady} derived a new model that neglects nonlinear
contributions in the moment equations, so that the additional evolution equations become much simpler. The model was thus called shallow water linearized moment
equations (SWLME). Thanks to its simplified structure, the SWLME allowed to study analytically
the steady states and equilibrium variables, which could be seen as an extension to the known SWE case. Those analytically computed steady states were then used to design the first
well-balanced schemes in the context of shallow water moment models \cite{koellermeier2022steady,caballero2025semi}.

The goal of this paper is twofold: On one hand, we aim at a general
well-balanced high order numerical schemes for complex models involving
non-conservative products, for which it is not possible to know their steady states a priori.
On the other hand, we aim to show how such proposed schemes can be efficiently applied to the different shallow water moment models mentioned above (including SWME, HSWME, and SWLME) and perform a direct comparison
of how perturbations evolve.

To achieve both, we use the global flux idea, first introduced in
\cite{gascon2001construction,caselles2009flux,donat2011hybrid} and further developed in
\cite{chertock2022well,mantri2024fully,KAZOLEA2025106646,barsukow2025structure,barsukow2025genuinely}.
The main idea of global flux methods is to build a quasi-conservative hyperbolic system starting
from a balance law by integrating the source term in a global definition of the flux.
In \cite{ciallella2023arbitrary} a promising strategy to construct high order
well-balanced finite volume schemes for the SWE was developed by combining the flux globalization approach with a
tailored quadrature based on the high order accurate WENO reconstruction.
In this work, we show how to extend the approach in \cite{ciallella2023arbitrary} to the
case of hyperbolic balance laws with non-conservative products.
This leads to a general scheme that allows to study the equilibria preservation of different shallow water moment models.

With a range of numerical experiments, the convergence, accuracy, and flexibility of the new schemes are demonstrated.
The results show that discretization errors for different shallow water moment models are comparable with the results obtained for the SWE in \cite{ciallella2023arbitrary},
which were several orders of magnitude smaller than those obtained with standard WENO schemes.
This holds true for higher order moment models as well, where our new scheme is able to capture the evolution
of small perturbations of the steady states with the same accuracy. A final comparison between
several moment models with both non-flat bathymetry and friction terms is presented to demonstrate the flexibility of our method.

The rest of the paper is organized as follows:
In \cref{se:SW_Equation}, we recall the standard SWE, their steady states, and their quasi-conservative global flux formulation.
In \cref{sec:moment}, we introduce the shallow water moment models, their steady states (when available), and their new global flux formulation.
In \cref{sec_Space}, we introduce the high order spatial discretization including the finite volume discretization, the global flux high order quadrature,
the WENO reconstruction, and the treatment of source terms and non-conservative products.
In \cref{sec:time}, we recall the time discretization used for the numerical simulations.
In \cref{sec:numerics}, we present the numerical tests to show the optimal convergence and flexibility
of our scheme when applied to different families of shallow water moment models.
Finally, in \cref{sec:conclusions}, we draw our conclusions and discuss future perspectives.

\section{Shallow Water Equations (SWE)}\label{se:SW_Equation}
\begin{figure}
	\centering
	\subfigure[Physical variables]{
	\includegraphics{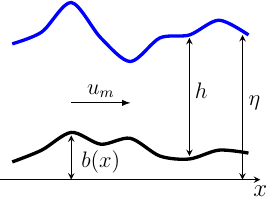}
	}\qquad\qquad
	\subfigure[Velocity]{
	\includegraphics{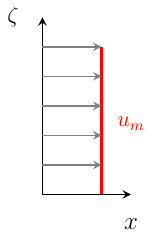}
	}
	\caption{Shallow water equations: model variables.}\label{fig:SWE}
\end{figure}

We start by describing the shallow water equations (SWE), which assume a constant velocity $u(x,\zeta,t)=u_m(x,t)$ over the scaled depth $\zeta \in [0,1]$ of the water column.
The simplified model resulting from this assumption can be recast as a system of hyperbolic balance laws as
\begin{equation}\label{eq:CL}
\partial_t \bu + \partial_x \bF(\bu) = \bS(\bu,x),
\end{equation}
with conserved variables $\bu$, flux $\bF(\bu)$, and source term $\bS(\bu,x)$ given by
\begin{equation}\label{eq:SWE}
	\bu=\begin{bmatrix} h \\ h u_m  \end{bmatrix},\quad
	\bF(\bu)=\begin{bmatrix} h u_m \\ h u_m^2 +g\frac{h^2}{2} \end{bmatrix},\quad
	\bS(\bu,x)=  -gh\begin{bmatrix} 0 \\ \partial_x b(x) \end{bmatrix}
	-\frac{\nu}{\lambda} \mathbf{P}(\bu),
\end{equation}
where $h$ represents the relative water height, $hu_m$ is the discharge where $u_m$ is the vertically averaged velocity,
$g$ is the gravitational acceleration, and $b(x)$ is the given local bathymetry, as shown in \cref{fig:SWE}.
The parameters $\nu$ and $\lambda$ represent the kinematic viscosity and slip length, respectively, for a Newtonian fluid with bottom slip law.
The right hand side friction term $\mathbf{P}(\bu)$ for the classical SWE simply reads
\begin{equation}\label{eq:frictionSW}
\mathbf{P}(\bu) = \begin{bmatrix} 0 \\ u_m \end{bmatrix}.
\end{equation}
It is also convenient to introduce the free surface water level $\eta:=h+b$.

In non-conservative form, the SWE can be written as
\begin{equation}\label{eq:nonCL}
	\partial_t \bu + \bA(\bu) \partial_x \bu = \bS(\bu,x),
\end{equation}
where $\bA(\bu)$ is the Jacobian matrix of the flux $\bF(\bu)$ with respect to the variables $\bu$:
\begin{equation}\label{eq:Jacobian}
	\bA(\bu) = \begin{bmatrix} 0 & 1 \\ -u_m^2 + gh & 2u_m \end{bmatrix}.
\end{equation}	
The eigenvalues of the Jacobian matrix are given by
\begin{equation}\label{eq:eigenvalues}
	\lambda_{1,2} = u_m \pm \sqrt{gh} .
\end{equation}
Following \cite{ciallella2023arbitrary}, the global flux method for the SWE
consists in writing \eqref{eq:CL} in the following equivalent quasi-conservative form:
\begin{equation}\label{eq:globalCL}
	\partial_t \bu + \partial_x \bG(\bu,x) = 0 \quad\text{ such that }\quad \bG(\bu,x)=\begin{bmatrix} h u_m \\ K_m\end{bmatrix}=\begin{bmatrix} h u_m \\ h u_m^2 +g\frac{h^2}{2}+\mathcal R_m \end{bmatrix},
\end{equation}
having set
\begin{equation}\label{eq:globalR}
	\mathcal R_m(x,t) := -\int_{x_0}^x S_2(\bu,\xi) \diff{\xi}= g \int_{x_0}^x \left[ h(\xi,t) \partial_\xi b(\xi) + \frac{\nu}{\lambda} u_m(\xi,t) \right]  \diff{\xi},
\end{equation}
where $x_0$ is a reference point, in general taken as the left boundary of the domain so that $\mathcal R_m$ is a local discrete approximation of the integral of the source term.

In this case both components of the global flux,  $hu_m$ and $K_m$, are steady discrete equilibrium variables since we can simply infer that
\begin{equation}\label{eq:steady}
	\partial_t \bu = 0 \quad\Longleftrightarrow\quad \bG(\bu,x) = \text{const.}
\end{equation}

In particular, the family of exact steady states for the SWE, with flat topography and no friction, fulfills
\begin{equation}
\begin{cases}
	\displaystyle \partial_x \left(hu_m\right) &= 0, \\
	\displaystyle h\partial_x \left(\frac{1}{2}u_m^2 + gh \right) &= 0,	
\end{cases}
\Longleftrightarrow
\begin{cases}
	 \displaystyle hu_m &= \text{const.}, \\
	 \displaystyle \frac{1}{2}u_m^2 + gh &= \text{const.}
\end{cases}
\end{equation}
Similarly, for a smooth frictionless flow with a non-zero bathymetry, the steady states fulfills
\begin{equation}
	\begin{cases}
		 \displaystyle hu_m &= \text{const.}, \\
		 \displaystyle \frac{1}{2}u_m^2 + g(h+b) &= \text{const.},	
	\end{cases}
\end{equation}
where the momentum equation has been manipulated, considering the constant discharge $hu_m$, as follows
%
\begin{align*}
0 &= \partial_x\left(hu_m^2 + g\frac{h^2}{2}\right)+gh\partial_x b = hu_m \partial_x u_m + gh \partial_x h + gh \partial_x b = h \partial_x \left(\frac12 u_m^2 + g(h+b)\right).
\end{align*}
%
%
\section{Shallow water moment equations (SWME)}\label{sec:moment}
The shallow water moment equations (SWME) generalize the SWE by permitting vertical variation of the water velocity profile \cite{kowalski2019moment}.
This is achieved by assuming a polynomial expansion of the velocity profile, which is then truncated to a finite number of terms $N$:
\begin{equation}\label{eq:momentvel}
u(x,\zeta,t) = u_m(x,t) + \sum_{k=1}^{N} \alpha_k(x,t) \phi_k(\zeta),
\end{equation}
where $u_m(x,t)$ is the mean horizontal velocity, used also in the standard SWE, $\alpha_k(x,t)$ are the polynomial coefficients or so-called \emph{moments} of the velocity profile and
$\phi_k(\zeta)$ are scaled, orthogonal Legendre basis functions, defined by
\begin{equation}
\phi_k(\zeta) = \frac{1}{j!}\frac{\diff{}^j}{\diff{\zeta}^j}\left(\zeta-\zeta^2\right)^j \quad \text{for } j = 1,\ldots,N.
\end{equation}
From \cref{eq:momentvel} we see that a larger \emph{order} $N$ leads to more allowed variations in the velocity profile, potentially increasing accuracy of free-surface simulations.
In this paper, we mainly consider the cases $N=1$, called the first-order or linear model, and $N=2$, called the second order or quadratic model. However, our method is applicable for general $N$.

\begin{figure}
	\centering
	\subfigure[]{
	\includegraphics{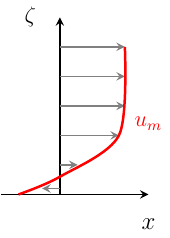}
	}\qquad\qquad
	\subfigure[]{
	\includegraphics{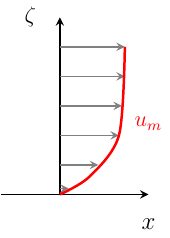}
	}
	\caption{Vertically varying velocity profiles with a change of sign (left), and without a change of sign (right). When a change of sign occurs, a vortex can form in the flow, which breaks the shallow water assumption.}\label{fig:vortexmoment}
\end{figure}

As an example, the first order system with $N=1$ leads to a velocity profile defined as a linear polynomial depending on mean velocity $u_m$ and first moment $\alpha_1$ as
\begin{equation}\label{eq:veldistribution}
u(x,\zeta,t) = u_m(x,t) + \left(1-2\zeta\right)\alpha_1(x,t), \quad\text{ with } \zeta = \frac{z-b}{h},
\end{equation}
where $z=h+b$ (or $\zeta = 1$) at the top, and $z=b$ (or $\zeta = 0$) at the bottom of the water column.
Hence, the velocity will take the following values at the top and bottom of the water column:
\begin{equation}
u(z=h+b) = u_m - \alpha_1 \quad \text{and} \quad u(z=b) = u_m + \alpha_1.
\end{equation}
It is important to notice that $u(\zeta)$ should have the same sign on the whole water column, otherwise the classical shallow water assumption is not valid anymore,
i.e.\ a vortex may form as shown in \cref{fig:vortexmoment}. Therefore, we require $|\alpha_1(x,t)| \leq |u_m(x,t)|$.

Evolution equations for the moments $\alpha_k$, for $k = 1,\ldots, N$ are derived by taking higher-order averages of the underlying incompressible Navier-Stokes equations. The resulting equations include non-conservative products $\bB(\bu)\partial_x \bu$ meaning that the SWME derived in \cite{kowalski2019moment} can be written in the following compact form
\begin{equation}\label{eq:CL-NC}
	\partial_t \bu + \partial_x \bF(\bu)  = \bB(\bu)\partial_x \bu + \bS(\bu,x).
\end{equation}

The dimensions and entries of the variables $\bu = (h, h u_m, h \alpha_1, \ldots, h \alpha_N)^T \in \mathbb{R}^{N+2}$, the flux $\bF \in \mathbb{R}^{N+2}$, the non-conservative products $\bB  \in \mathbb{R}^{(N+2)\times (N+2)}$ and the source term $\bS  \in \mathbb{R}^{N+2}$ depend on the truncation order $N$. We note that the SWME for different $N$ are well-studied and can be found in \cite{kowalski2019moment,koellermeier2020analysis,koellermeier2022steady}. Below we will give the explicit formulas for the low-order models with $N=1$ and $N=2$.

\subsection{First order shallow water moment equations (SWME1)}
When considering a linear velocity profile, i.e.\ $N=1$, the first order shallow water moment equations SWME1 from \cite{kowalski2019moment} read
as in \cref{eq:CL-NC} where $\bu$, $\bF, \bS$ and $\mathbf{P}$ are defined as follows
\begin{align}\label{eq:SWME1}
	\bu=\begin{bmatrix} h \\ h u_m \\ h \alpha_1  \end{bmatrix},\quad
	\bF(\bu)=\begin{bmatrix} h u_m \\ h u_m^2 +g\frac{h^2}{2} + \frac13 h \alpha_1^2 \\ 2 h u_m \alpha_1 \end{bmatrix},
\end{align}
\begin{align}
\bS(\bu,x)=  -gh\begin{bmatrix} 0 \\ 0 \\ \partial_x b(x) \end{bmatrix}
	-\frac{\nu}{\lambda} \mathbf{P}(\bu), \quad
	\mathbf{P}(\bu) = \begin{bmatrix} 0 \\ u_m + \alpha_1 \\ 3 \left(u_m + \alpha_1 + 4\frac{\lambda}{h}\alpha_1 \right) \end{bmatrix},
\end{align}
while $\bB$ represents the matrix of non-conservative products given by
\begin{align}\label{eq:SWME1-B}
	\bB(\bu)= \begin{bmatrix} 0 & 0 & 0 \\ 0 & 0 & 0 \\ 0 & 0 & u_m  \end{bmatrix}.
\end{align}
The flux Jacobian can be directly computed and reads
\begin{align*}
	\partial_{\bu} \bF = \begin{bmatrix} 0 & 1 & 0 \\ -u_m^2+gh-\frac{\alpha_1^2}{3} & 2 u_m & \frac{2\alpha_1}{3} \\ -2 u_m \alpha_1 & 2 \alpha_1 & 2 u_m  \end{bmatrix},
\end{align*}
and the system matrix, which allows us to rewrite the system as \eqref{eq:nonCL}, is
\begin{align*}
	\bA = \partial_{\bu} \bF - \bB = \begin{bmatrix} 0 & 1 & 0 \\ -u_m^2+gh-\frac{\alpha_1^2}{3} & 2 u_m & \frac{2\alpha_1}{3} \\ -2 u_m \alpha_1 & 2 \alpha_1 & u_m  \end{bmatrix}.
\end{align*}
The SWME1 are strictly hyperbolic and have three distinct eigenvalues
$$ \lambda_{1,2} = u_m \pm \sqrt{gh + \alpha_1^2} \quad \text{and} \quad \lambda_3 = u_m .$$

The first study regarding steady states of the SWME1 can be found in \cite{koellermeier2022steady}, where a closed form for the general equilibria was provided.
In particular, we recall that, for flat topography and no friction, the steady states of the SWME1 fulfill
\begin{equation}
	\begin{cases}
		\displaystyle \partial_x \left(hu_m\right) &= 0, \\
		\displaystyle \partial_x \left(hu_m^2 + g\frac{h^2}{2} + \frac13 h \alpha^2_1 \right) &= 0,\\
		\displaystyle \partial_x \left(2 h u_m \alpha_1 \right) &= u_m \partial_x(h\alpha_1).
	\end{cases}
\end{equation}
When considering the constant discharge $hu_m$, the last equation can be recast as
%
\begin{align*}
0 &= \partial_x \left(2 h u_m \alpha_1 \right) - u_m \partial_x(h\alpha_1) = \frac{u_m}{h} \partial_x (\alpha_1) - \frac{u_m \alpha_1}{h^2} \partial_x h = u_m \partial_x \left( \frac{\alpha_1}{h} \right),
\end{align*}
which is fulfilled when
$$ u_m = 0 \quad \text{ or } \quad  \frac{\alpha_1}{h} = \text{const.} $$

With the non-trivial conditions coming from the first and third equation ($hu_m=\text{const.}$ and $\alpha_1/h=\text{const.}$),
we can recast the momentum equation as follows
%
\begin{align*}
  0 &= \partial_x \left(hu_m^2 + g\frac{h^2}{2} + \frac13 h \alpha^2_1 \right) = h \partial_x \left( \frac12 u_m^2 + g h + \frac12 \alpha_1^2 \right).
\end{align*}
Similarly, if the bathymetry term is present in the equation, the non trivial steady states are described by
\begin{equation}\label{eq:SWME1-steady}
	\begin{split}
	h u_m &= \text{const.}  \\
	\frac12 u_m^2 + g (h+b) + \frac12 \alpha_1^2 &= \text{const.}  \\
	\frac{\alpha_1}{h} &= \text{const.}
	\end{split}
\end{equation}
From the first and third equation of \eqref{eq:SWME1-steady}, we can introduce the following constants
$$ h u_m = C_0 \qquad  \text{and} \qquad \frac{\alpha_1}{h} = C_1. $$
Then, the second equation of \eqref{eq:SWME1-steady} can be recast as
\begin{equation}\label{eq:SWME1-steadyH}
	\frac{C_0^2}{2 g h^2} + h + b + \frac{C_1^2 h^2}{2 g} = \text{const.}
\end{equation}
Thanks to the closed form of the equilibria, we are able to compute the exact value of $h$ by fixing the right-hand side
of \cref{eq:SWME1-steadyH}, for instance using the data at the left boundary of the domain similarly to \cite{delestre2013swashes}.
This directly allows us to compute $h$ by solving the following fourth order nonlinear equation:
\begin{equation}\label{eq:SWME1-exactH}
 \frac{C_1^2}{2 g} h^4 + h^3 + \left( b - \frac{C_0^2}{2 g h_L^2} - h_L - b_L - \frac{C_1^2 h_L^2}{2 g} \right) h^2 + \frac{C_0^2}{2 g} = 0 .
\end{equation}

For the SWME1 model, the global flux formulation of the system reads
\begin{equation}\label{eq:globalSWME1}
	\partial_t \bu + \partial_x \bG(\bu,x) = 0 \quad\text{ such that }\quad \bG(\bu,x)=\begin{bmatrix} h u_m \\ K_m \\ K_1 \end{bmatrix}=\begin{bmatrix} h u_m \\ h u_m^2 +g\frac{h^2}{2} + \frac13 h \alpha_1^2 +\mathcal R_m \\ 2 h u_m \alpha_1 + \mathcal R_1 \end{bmatrix},
\end{equation}
with $\bR = \begin{bmatrix} \mathcal R_m \\ \mathcal R_1 \end{bmatrix}$:
\begin{equation}\label{eq:globalR_SWME1}
\begin{split}
	\mathcal R_m(x,t) &:= \int_{x_0}^x \left[ g h \partial_\xi b + \frac{\nu}{\lambda} (u_m + \alpha_1) \right]  \diff{\xi}, \\
	\mathcal R_1(x,t) &:= \int_{x_0}^x \left[ -u_m \partial_\xi (h \alpha_1) + 3\frac{\nu}{\lambda} \left(u_m + \left(1+4\frac{\lambda}{h}\right)\alpha_1 \right) \right]  \diff{\xi},
\end{split}
\end{equation}
where, for convenience, we dropped the $(\xi,t)$ dependence.
It can be noticed in \cref{eq:globalR_SWME1} that the non-conservative product is treated as an additional source term in the global flux formulation.

\subsection{Second order shallow water moment equations (SWME2)}
The second order shallow water moment equations (SWME2) can be derived from the hypothesis that the velocity profile is a polynomial of degree two along the vertical direction, i.e., $N=2$.
The SWME2 can be recast in compact form \eqref{eq:CL-NC} with the following definitions

\begin{align}\label{eq:SWME2}
	\bu=\begin{bmatrix} h \\ h u_m \\ h \alpha_1 \\ h \alpha_2  \end{bmatrix},\quad
	\bF(\bu)=\begin{bmatrix} h u_m \\ h u_m^2 +g\frac{h^2}{2} + \frac13 h \alpha_1^2 + \frac15 h \alpha_2^2 \\ 2 h u_m \alpha_1 + \frac45 h \alpha_1 \alpha_2 \\ 2 h u_m \alpha_2 + \frac23 h \alpha_1^2 + \frac27 h \alpha_2^2  \end{bmatrix},\quad
	\mathbf{P}(\bu) = \begin{bmatrix} 0 \\ u_m + \alpha_1 + \alpha_2 \\ 3 \left(u_m + \alpha_1 + \alpha_2 + 4\frac{\lambda}{h}\alpha_1 \right) \\ 5 \left(u_m + \alpha_1 + \alpha_2 + 12\frac{\lambda}{h}\alpha_2 \right) \end{bmatrix},
\end{align}
with the following matrix $\bB$ of non-conservative products
\begin{align}\label{eq:SWME2-B}
	\bB(\bu)= \begin{bmatrix} 0 & 0 & 0 & 0 \\ 0 & 0 & 0 & 0 \\ 0 & 0 & u_m - \frac{\alpha_2}{5} & \frac{\alpha_1}{5} \\ 0 & 0 & \alpha_1 & u_m + \frac{\alpha_2}{7} \end{bmatrix}.
\end{align}
The flux Jacobian, in this case, is
\begin{align*}
	\partial_{\bu} \bF = \begin{bmatrix} 0 & 1 & 0  & 0 \\ -u_m^2+gh-\frac{\alpha_1^2}{3}-\frac{\alpha_2^2}{5} & 2 u_m & \frac{2\alpha_1}{3} & \frac{2\alpha_2}{5} \\
		-2 u_m \alpha_1 - \frac45 \alpha_1 \alpha_2 & 2 \alpha_1 & 2 u_m + \frac{4\alpha_2}{5} & \frac{4\alpha_1}{5} \\
	    -2 u_m \alpha_2 - \frac23 \alpha^2_1 - \frac27 \alpha^2_2 & 2 \alpha_2 & \frac{4\alpha_1}{3} & 2 u_m + \frac{4\alpha_2}{7} \end{bmatrix},
\end{align*}
and the system matrix reads
\begin{align*}
	\bA = \partial_{\bu} \bF - \bB = \begin{bmatrix}  0 & 1 & 0  & 0 \\ -u_m^2+gh-\frac{\alpha_1^2}{3}-\frac{\alpha_2^2}{5} & 2 u_m & \frac{2\alpha_1}{3} & \frac{2\alpha_2}{5} \\
		-2 u_m \alpha_1 - \frac45 \alpha_1 \alpha_2 & 2 \alpha_1 & u_m + \alpha_2 & \frac{3\alpha_1}{5} \\
	    -2 u_m \alpha_2 - \frac23 \alpha^2_1 - \frac27 \alpha^2_2 & 2 \alpha_2 & \frac{\alpha_1}{3} & u_m + \frac{3\alpha_2}{7} \end{bmatrix}.
\end{align*}

\begin{figure}
	\centering
	\includegraphics[width=0.3\textwidth]{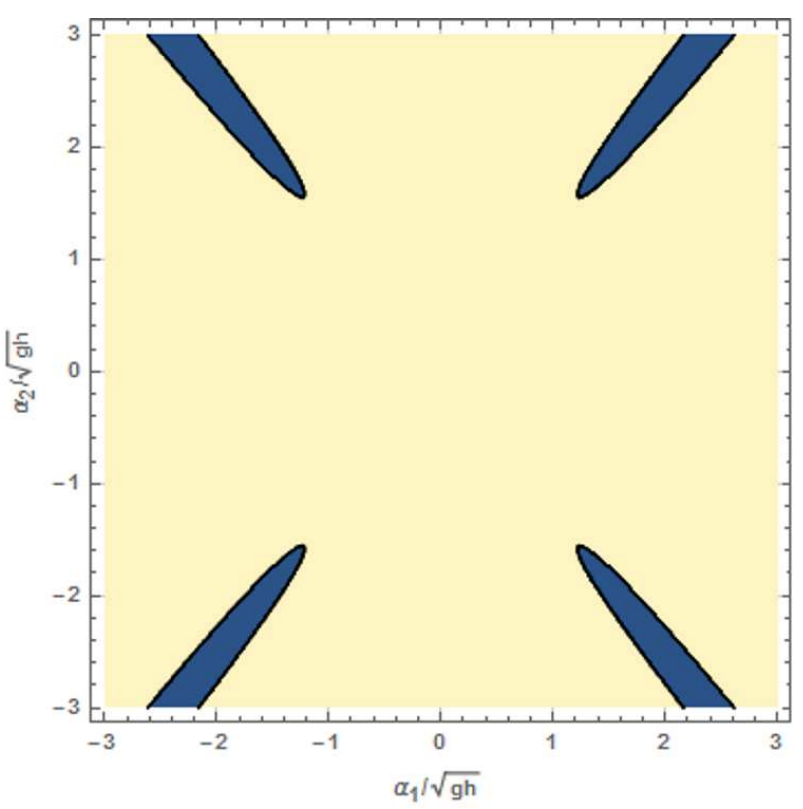}
	\caption{Second order shallow water moment model: loss of hyperbolicity (blue) and hyperbolic region (beige) taken from \cite{koellermeier2020analysis}.}\label{fig:hyp-SWME2}
\end{figure}

Hyperbolicity, i.e., the existence of real eigenvalues and a full set of eigenvectors of the matrix $\bA$, is an important requirement in the analysis and simulation of such models.
Already with this small number of moments $N=2$, the SWME2 may loose hyperbolicity, as analyzed in \cite{koellermeier2020analysis}.
In particular, it was shown that the model is hyperbolic only for certain states depending on the coefficients $\alpha_1$ and $\alpha_2$, which should not be very large, see \cref{fig:hyp-SWME2}.

However, it should be noticed that a loss of hyperbolicity is not uncommon when considering more complicated physical models. For instance, the two-layer shallow water model \cite{lawrence1993hydraulics,castro2004numerical} also experience
hyperbolicity loss due to the development of shear instabilities.
Beyond the field of shallow water flows, loss of hyperbolicity is also encountered in moment models of rarefied gases \cite{Koellermeier2014}, which has led to regularized alternatives including non-conservative products \cite{Koellermeier2017a}.

In a similar fashion hyperbolic regularizations of shallow water moment models have been proposed in the literature to obtain more stable results while achieving similar accuracy as the original SWME models for $N\geq 2$, see \cite{koellermeier2020analysis,koellermeier2022steady}. We will outline two of those in the next section.

Regardless of hyperbolicity, the global flux formulation of the SWME2 reads
\begin{equation}\label{eq:globalSWME2}
	\bG(\bu,x)=\begin{bmatrix} h u_m \\ K_m \\ K_1 \\ K_2 \end{bmatrix}=\begin{bmatrix} h u_m \\ h u_m^2 +g\frac{h^2}{2} + \frac13 h \alpha_1^2 + \frac15 h \alpha_2^2 +\mathcal R_m \\
		2 h u_m \alpha_1 + \frac45 h \alpha_1 \alpha_2 + \mathcal R_1 \\ 2 h u_m \alpha_2 + \frac23 h \alpha_1^2 + \frac27 h \alpha^2_2 + \mathcal R_2 \end{bmatrix},
\end{equation}
where $\bR = \begin{bmatrix} \mathcal R_m \\ \mathcal R_1  \\ \mathcal R_2 \end{bmatrix}$ are given by
\begin{equation}\label{eq:globalR_SWME2}
	\begin{split}
		\mathcal R_m(x,t) &:= \int_{x_0}^x \left[ g h \partial_\xi b + \frac{\nu}{\lambda} P_2 \right]  \diff{\xi}, \\
		\mathcal R_1(x,t) &:= \int_{x_0}^x \left[ -B_{33} \partial_\xi (h \alpha_1) -B_{34} \partial_\xi (h \alpha_2) + \frac{\nu}{\lambda} P_3 \right]  \diff{\xi}, \\
		\mathcal R_2(x,t) &:= \int_{x_0}^x \left[ -B_{43} \partial_\xi (h \alpha_1) -B_{44} \partial_\xi (h \alpha_2) + \frac{\nu}{\lambda} P_4 \right]  \diff{\xi},
	\end{split}
\end{equation}
where we used the following definitions
$$ \bB =\begin{bmatrix} 0 & 0 & 0 & 0 \\ 0 & 0 & 0 & 0 \\ 0 & 0 & B_{33} & B_{34} \\ 0 & 0 & B_{43} & B_{44}\end{bmatrix} ,\qquad \mathbf{P}(\bu) = \begin{bmatrix} P_1 \\ P_2 \\ P_3 \\ P_4  \end{bmatrix}. $$
We will also use the same notation below to have more compact formulas.

\subsection{Second order hyperbolic shallow water moment equations (HSWME2)}

The breakdown of hyperbolicity, and the possible source of instabilities attributed to that, motivated the hyperbolic regularization in \cite{koellermeier2020analysis}, which changes the system matrix $\bA$ such that the system is hyperbolic.
In the general case of $N\geq 1$, one of the studied approaches is to set to zero in $\bA$ all
coefficients $\alpha_i$, for $i>1$. This provides a new family of hyperbolic models,
called the hyperbolic shallow water moment equations (HSWME).

For $N=2$, the system matrix for the second order hyperbolic shallow water moment equations (HSWME2) reads
$$
\bA =  \begin{bmatrix}  0 & 1 & 0  & 0 \\
	-u_m^2+gh-\frac{\alpha_1^2}{3} & 2 u_m & \frac{2\alpha_1}{3} & 0 \\
	-2 u_m \alpha_1  & 2 \alpha_1 & u_m  & \frac{3\alpha_1}{5} \\
	-\frac23 \alpha_1^2  & 0 & \frac{\alpha_1}{3} & u_m  \end{bmatrix}.
$$
In this case, the real propagation speeds can be explicitly computed as
\begin{equation}\label{eq:HSWME2speeds}
\lambda_{1,2} = u_m \pm \sqrt{gh + \alpha_1^2}\quad\text{ and } \quad \lambda_{3,4} = u_m \pm \sqrt{\frac15 \alpha_1^2}.
\end{equation}
Although the conservative variables $\bu$ and the source term $\bS$ stay the same as in \cref{eq:SWME2},
the conservative flux $\bF$ and the matrix of non-conservative product $\bB$ take a different shape:
\begin{align}\label{eq:HSWME2}
	\bF(\bu)=\begin{bmatrix} h u_m \\ h u_m^2 +g\frac{h^2}{2} + \frac13 h \alpha_1^2  \\ 2 h u_m \alpha_1 \\ \frac23 h \alpha_1^2   \end{bmatrix},\quad
	\bB(\bu)= \begin{bmatrix} 0 & 0 & 0 & 0 \\ 0 & 0 & 0 & 0 \\ 0 & 0 & u_m & -\frac{3\alpha_1}{5} \\ 0 & 0 & \alpha_1 & -u_m \end{bmatrix}.
\end{align}
For this model, the global flux formulation reads
\begin{equation}\label{eq:globalHSWME2}
	\bG(\bu,x)=\begin{bmatrix} h u_m \\ K_m \\ K_1 \\ K_2 \end{bmatrix}=\begin{bmatrix} h u_m \\ h u_m^2 +g\frac{h^2}{2} + \frac13 h \alpha_1^2 +\mathcal R_m \\
		2 h u_m \alpha_1 + \mathcal R_1 \\ \frac23 h \alpha_1^2 + \mathcal R_2 \end{bmatrix},
\end{equation}
where the definitions of $\mathcal R_m, \mathcal R_1, \mathcal R_2$ are the same as in \cref{eq:globalR_SWME2}
with the appropriate coefficients for the non-conservative products.

As studied in \cite{koellermeier2022steady}, for complex models such as SWME2 and HSWME2, it is not possible to write steady states
in a closed form as done for the SWME1 model in \cref{eq:SWME1-steady}.
This lack of an analytical formula for the steady states is mainly related to the fact that the number of non-conservative products is large, which make them hard to manipulate.
This might pose several problems in the design of well-balanced numerical schemes as discussed previously in \cite{koellermeier2020analysis},
and then in \cite{koellermeier2022steady}.
For this reason \cite{koellermeier2022steady} developed a new class of hyperbolic models,
called the shallow water linearized moment equations (SWLME).
With the goal of developing well-balanced schemes, these models have a simpler structure that allows for an analytical computation of steady states in closed form.

\subsection{Second order shallow water linearized moment equations (SWLME2)}

As described in \cite{koellermeier2022steady}, the derivation of the shallow water linearized moment (SWLME) models originates from the goal of having a closed form of non-trivial steady states,
which is then employed to develop a numerical scheme that preserves those steady states.
In particular, this can be managed as long as there are not that many non-conservative terms in the model.

In \cite{koellermeier2022steady} this is achieved by assuming small expansion coefficients $\alpha_i = \mathcal{O}(\epsilon) $ and then neglecting higher-order contributions $\mathcal{O}(\epsilon^2)$ to simplify the higher-order moment equations.
This simplification leads to the new hyperbolic SWLME, which allow to recover a closed form of the steady states also for arbitrary $N$.

As we recall below, it results in right-hand side non-conservative terms that do no longer contain coupling between the $\alpha_i$,
basically transforming high order models to a simple extensions of the first order model, for which the steady states are known.

The general SWLME with arbitrary number of moments $N$ from \cite{koellermeier2022steady} can then be written in general form. Here we only consider the model for $N=2$, where the variables $\bu$ and the source term $\bS$ are the same as in \cref{eq:SWME2},
while $\bF$ and $\bB$ read
\begin{align}\label{eq:SWLME2}
	\bF(\bu)=\begin{bmatrix} h u_m \\ h u_m^2 +g\frac{h^2}{2} + \frac13 h \alpha_1^2 + \frac{1}{5} h \alpha_2^2\\ 2 h u_m \alpha_1 \\ 2 h u_m \alpha_2 
\end{bmatrix},\quad
	\bB(\bu)= \text{diag}(0,\,0,\,u_m,\,u_m).
\end{align}
In this case, the system matrix is
$$
\bA =  \begin{bmatrix}  0 & 1 & 0  & 0 \\
	-u_m^2+gh-\frac{\alpha_1^2}{3}-\frac{\alpha_2^2}{5} & 2 u_m & \frac{2\alpha_1}{3}  & \frac{2\alpha_2}{5} \\
	-2 u_m \alpha_1  & 2 \alpha_1 & u_m  &  0    \\
	-2 u_m \alpha_2  & 2 \alpha_2 &  0    & u_m \end{bmatrix}.
$$
For the SWLME2 model, the real propagation speeds can be explicitly computed as
\begin{equation}\label{eq:SWLME2speeds}
\lambda_{1,2} = u_m \pm \sqrt{gh + \alpha_1^2 + \frac{3}{5} \alpha_2^2}\quad\text{ and } \quad \lambda_{3,4} = u_m.
\end{equation}
Thus, the system is hyperbolic for all values of the coefficients $\alpha_1, \alpha_2$, and for positive water heights $h$.
Following the same reasoning of the SWME1 model, it is possible to retrieve the moving equilibria of the SWLME2 as
\begin{equation}\label{eq:steady_SWLME}
	\begin{split}
		h u_m &= \text{const.},  \\
		\frac12 u_m^2 + g (h+b) + \frac{1}{2}\alpha_1^2 + \frac{3}{10} \alpha_2^2 &= \text{const.},  \\
		\frac{\alpha_i}{h} &= \text{const.},\qquad \text{for } i=1,2.	
	\end{split}
\end{equation}

The global flux formulation for this model reads
\begin{equation}\label{eq:globalSWLME2}
	\bG(\bu,x)=\begin{bmatrix} h u_m \\ K_m \\ K_1 \\ K_2 \end{bmatrix}=\begin{bmatrix} h u_m \\ h u_m^2 +g\frac{h^2}{2} + \frac13 h \alpha_1^2 + \frac{1}{5} h \alpha_2^2 + \mathcal R_m \\
		2 h u_m \alpha_1 + \mathcal R_1 \\ 2 h u_m \alpha_2 + \mathcal R_2 \end{bmatrix},
\end{equation}
where $\mathcal R_m$ is the same as in \cref{eq:globalR_SWME2} and the $\mathcal R_i$ are defined as
\begin{equation}\label{eq:globalR_SWLME}
	\mathcal R_i(x,t) := \int_{x_0}^x \left[ - u_m \partial_\xi (h \alpha_i) + \frac{\nu}{\lambda} P_{i+2} \right]  \diff{\xi},\quad i=1,2.
\end{equation}

In this work, which is focused on the numerical method for general models, we only write the SWLME for $N=2$. However, the definition of the SWLME for general $N$ including their explicit eigenvalues and steady states can be found in \cite{koellermeier2022steady}. We note that all methods of this paper can readily be extended for larger $N$ using the respective higher order models SWME, HSWME, SWLME.

\section{Space discretization: Global Flux Finite Volume method}\label{sec_Space}

In this section, we present the space discretization of the system of hyperbolic balance laws \eqref{eq:CL-NC} with global fluxes.
The goal is to further develop the high order WENO method presented in the context of
the classical shallow water equations in \cite{ciallella2023arbitrary}, and extend it to more complex models with larger number of equations and involving
non-conservative products. Although our new approach is extremely general and can be applied to obtain
high order finite volume well-balanced methods for all kinds of hyperbolic systems with non-conservative products,
we will specifically focus on the construction of equilibria preserving methods for complex shallow water moment models in this paper.
The flexibility of the method is highlighted by the fact that it can be applied straightforwardly to all the shallow water
moment models, even those involving a strong nonlinear structure that does not allow us to compute the steady states in closed form.

The hyperbolic system considered herein is solved by means of the method of lines, hence space and time can be treated independently.
The computational domain $\boldsymbol{\Omega}$ is discretized into $N_x$ control volumes $\Omega_i = [\xin,\xip]$ of size $\Delta x$
centered at $x_i=i\Delta x$ with $i=i_\ell,\ldots,i_r$.

For the control volume $\Omega_i$ we can define the cell average at time $t$:
\begin{equation}
	\bar\bU_i(t):=\frac{1}{\Delta x}\int_{\xin}^{\xip} \bu(x,t) \diff{x}.
\end{equation}
The semi-discrete finite volume scheme for the system~\eqref{eq:globalCL} reads
\begin{equation}\label{eq:FV}
	\frac{\diff{\bar\bU_i}}{\diff{t}} + \frac{1}{\Delta x}(\widehat \bG_{\iip}-\widehat \bG_{\iin}) = 0 ,
\end{equation}
where $\widehat \bG_{\iip}$ is a numerical flux consistent with the global flux $\bG$.
\begin{remark}[Numerical global flux]
	It should be noticed that in order to achieve equilibria preservation,
	the numerical flux $\widehat \bG_{\iip}$ must only depend on the global flux $\bG$.
	This concept is crucial and the main difference from classical finite volume methods,
	where the numerical flux is defined as a central flux term plus a dissipation depending on conservative variables.
	This is due to the fact that, at equilibria, the global flux $\bG$ is constant while the conservative variables may vary.
	Hence, this definition provides preservation of steady states since the dissipation term vanishes at equilibria.
    \label{remark:gf}
\end{remark}
In this work, we are focusing on two different numerical global fluxes which will be compared by discussing both
the methodological point of view, and their impact on the numerical simulations.

The first one is a simple upwind flux, also used in \cite{ciallella2023arbitrary}, defined as
\begin{align}\label{eq:upwind flux}
	\widehat \bG^u_{\iip} &= (L^{-1}\Lambda^+L)_{\iip} \bG^L_{\iip}  + (L^{-1}\Lambda^-L)_{\iip} \bG^R_{\iip},
\end{align}
where, $\bG^{L,R}_{\iip}$ are the discontinuous reconstructed point values of the global flux
$\bG(\bU)$ respectively at the left and right side of the cell interface $x_{\iip}$.
$L$ is the matrix of the left eigenvectors computed from the system matrix $\bA$ of the hyperbolic
problem in the averaged state. $\Lambda^\pm$ correspond to the upwinding eigenvalues $\Lambda^\pm = \frac{\Lambda \pm \left| \Lambda \right|}{2}$.
For instance, in the case of the SWE, we would have
\begin{equation}
	\bA(\bU^*) = \begin{bmatrix}
		0 & 1\\ -u_*^2 +gh_* & 2u_*
	\end{bmatrix}, \qquad \text{with }\begin{cases}
		h_*=\frac{h^L+h^R}{2},\\
		u_*=\frac{u^L+u^R}{2}.
	\end{cases}
\end{equation}
Since the solution state used to compute the system matrix $\bA$ does not affect the well-balanced property of the scheme (see remark \ref{remark:gf}),
we can simply use the primitive reconstructed values at the interface.
Of course, also other types of averages can be used to compute the system matrix $\bA$, such as the Roe's average \cite{roe1981approximate}.
Other numerical global fluxes can be used as long as the numerical dissipation term depending on conservative
variables vanishes when the steady state is reached (see remark \ref{remark:gf}).
For example, this can be achieved using cutoff functions~\cite{chertock2018well},
however, some parameters must be carefully tuned to achieve equilibria preservation in that case.

The second numerical global flux is based on the idea of having a central part and a dissipative term,
similarly to what was done in \cite{barsukow2025genuinely}, where the latter is also defined as a function of the system matrix and global fluxes:
\begin{align}\label{eq:central flux}
	\widehat \bG^c_{\iip} &= \frac12 (\bG^L_{\iip}  +  \bG^R_{\iip}) - \frac{1}{|\lambda_{\text{max}}|}\bA(\bU^*) (\bG^R_{\iip}-\bG^L_{\iip}),
\end{align}
with $|\lambda_{\text{max}}|$ the spectral radius of $\bA$.
The main advantage of the latter is that it is easier to implement and less expensive than the upwind flux,
since it does not require the computation of the left and right eigenvectors.

To obtain the high order reconstructed values at the left and right side of interfaces,
we will use a high order WENO reconstruction technique on the cell averages of the global fluxes
coupled with a sophisticated high order treatment of both non-conservative products and source terms.
\subsection{Global flux high order quadrature}
In this section, we will discuss the high order quadrature of the global fluxes based on the WENO reconstruction.
In particular, to obtain a high order reconstruction at the interfaces of $\bG$, we need to compute the cell averages $\bar\bG_i$.
A consistent way to define them is to start from the cell averages of the conservative fluxes and the integral of non-conservative products and source terms:
\begin{equation}\label{eq:cellave global flux}
\bar\bG_i(\bu,x)=\bar\bF_i(\bu) + \bar\bR_i(\bu,x),
\end{equation}
with the cell averages defined using standard high order quadrature formulas (Gauss-Legendre in our case)
\begin{equation}\label{eq:flux_average}
	\bar\bF_i(\bu) = \sum_{q} w_q \bF(\tilde{\bu}(x_{i,q})) \quad\text{ and }\quad \bar\bR_i(\bu,x) = \sum_{q} w_q \bR(\tilde{\bu}(x_{i,q}),x_{i,q}),
\end{equation}
where $\{x_{i,q},w_q\}$ are high order quadrature points and weights used in the cell $[x_{i-1/2},x_{i+1/2}]$.
The notation $\tilde{\bu}(x_{i,q})$ indicates the high order WENO reconstruction of the conservative variables at the quadrature points $x_{i,q}$.
Below, we will discuss the definition of the integral term $\bR$ and how to compute it.

We can generally define the right-hand side of \cref{eq:CL-NC} as
$$\bH(\bU) = \bB(\bU)\partial_x \bU + \bS(\bU).$$
Here we provide the general formulation for hyperbolic systems with non-conservative products,
while an example of its application to shallow water moment models with source terms is given in the following sections.
Hence, to obtain the values at the quadrature points, we use a piecewise polynomial reconstruction of
the integral term $\bR$, and define it in a general recursive way:
\begin{equation}\label{eq: Rquad}
	\begin{split}
		\bR_{i,q} &= \bR^R_{\iin} -  \sum_\theta \left(\int_{\xin^R}^{x_{i,q}} L_\theta(x) \diff{x}\right)  \bH(\tilde{\bu}(x_{i,\theta})), \\
		&= \bR^R_{\iin} -  \sum_\theta \underbrace{\left(\int_{\xin^R}^{x_{i,q}} L_\theta(x) \diff{x}\right)}_{r_\theta^q}  \left(\bB(\tilde{\bu}(x_{i,\theta})) \sum_s L'_s(x_{i,\theta}) \tilde{\bu}(x_{i,s}) +  \bS(\tilde{\bu}(x_{i,\theta}))\right),
	\end{split}
\end{equation}
where $L_\theta$ are the Lagrangian polynomials associated to  the quadrature points, and their integrals $r_\theta^q$ are computed exactly.
Notice that $x_{i,q}$, $x_{i,\theta}$ and $x_{i,s}$ are the quadrature nodes on the interval $[x_{i-1/2},x_{i+1/2}]$.
As one may notice, the non-conservative products are also treated supposing a Lagrangian interpolation of the conservative variables:
$$ \partial_x\bu(x_{i,q}) = \sum_s L_s'(x_{i,q}) \tilde{\bu}(x_{i,s}). $$
Moreover, due to the finite volume formulation, we may have to deal with jumps of $\bR$ at the interfaces
and, for this reason, we introduce the notation $\bR^R_{\iip}$ and $\bR^L_{\iip}$ to indicate the values of $\bR$ at the right and left side of the interface $x_{\iip}$, respectively.

So far the notation is general and can be applied to any system of hyperbolic balance laws written in global flux form.
A special treatment of the bathymetry source term is needed to ensure exact well-balancing of the scheme for non-moving equilibria, e.g., lake at rest,
and will be further discussed in \cref{sec:larSourceDiscretization}.
Similarly, also the non-conservative product needs to be treated carefully to ensure consistency,
and will be detailed in \cref{sec:nonconservative}.

\begin{remark}[Boundary treatment of $\bR$]
The definition of $\bR$ requires an appropriate initial value $ \bR^R_{\iin} $.
When the source term is not acting at the left boundary of the computational domain, we can set it directly to zero.
However, when $\bR$ is already non-zero at the left boundary, we starts the integration from the ghost cells used to define the WENO polynomials at boundaries.
Since the ghost cells include the boundary conditions, this is enough to ensure that the integral $\bR$ is well defined at boundaries.
\end{remark}

To complete the definition of the iterative procedure, we need to link the values of $\bR$ at the left and right interfaces.
In particular, we set
\begin{equation}\label{eq:jumpR}
\bR^R_{\iip} = \bR^L_{\iip} + [\![ \bR ]\!]_{\iip},
\end{equation}
where we note that $\bR^L_{\iip}$ is
\begin{equation}\label{eq: Rinterface}
\begin{split}
	\bR^L_{\iip} \! &= \bR^R_{\iin} - \int_{\xin^R}^{\xip^L} \bH(\tilde\bU) \diff{x} \\
&=\bR^R_{\iin} -  \sum_\theta \left(\int_{\xin^R}^{\xip^L} L_\theta(x) \diff{x}\right)\!\!\!  \left(\bB(\tilde{\bu}(x_{i,\theta}))	 \sum_s L'_s(x_{i,\theta}) \tilde{\bu}(x_{i,s}) +  \bS(\tilde{\bu}(x_{i,\theta}))\right)\\
&=  \bR_{\iin}^R - \Delta x \bar{\bH}_i,
\end{split}
\end{equation}
where the last equality is obtained by the definition of the cell average.
To complete the algorithm, we need to provide a precise definition of the jumps $[\![ \bR ]\!]_{\iip}$.

\subsection{Weighted Essentially Non-Oscillatory (WENO) reconstruction}\label{sec:WENO}

In this section, we briefly recall the main ingredients of the WENO reconstruction \cite{jiang1996efficient,balsara2000monotonicity} used in this work.
In particular, we take advantage of the WENO method to perform high order reconstructions of the solution at quadrature points, and of global fluxes at the cell interfaces, starting from cell averages.

We briefly recall the  basics to obtain the polynomial of a scalar quantity $u_i(x)$ in cell $\Omega_i$.
Herein, we consider  polynomial reconstructions of order $p$, with $p$ odd and, to construct them, we select a high order stencil $\mathcal S_i$
of  $p$ cells centered in cell $\Omega_i$:
\begin{equation}
\mathcal S_i = \lbrace\Omega_j,\quad j = i-r+1, \dots, i+r-1 \rbrace,
\end{equation}
where $2r-1=p$. On each of these stencils, one constructs a high order  polynomial $P^{HO}$  fulfilling  the constraints
\begin{equation}
\frac{1}{\Delta x}\int_{x_{i-j-1/2}}^{x_{i-j+1/2}}P^{HO}(x) \,\diff{x} = u_{i-j}, \qquad j =-r+1,\dots, r-1,
\end{equation}
and $r$ low order polynomials $P^{LO}_m(x)$, $m=0,\dots ,r-1$, that fulfill
\begin{equation}
\frac{1}{\Delta x}\int_{x_{i-r+j+m-1/2}}^{x_{i-r+j+m+1/2}} P^{LO}_m(x)\diff{x} = u_{i-j+m}, \qquad j=1,\dots,r .
\end{equation}
The WENO method combines the low order polynomials to obtain a high order reconstruction, when the low order polynomials are non-oscillatory,
while it will weight more the least oscillatory polynomial in case some of these present oscillations.
For instance, the WENO reconstruction of order 5 (WENO5) uses a stencil of $p=5$, with $r=3$ low order reconstructions. \\
To achieve optimal accuracy, the linear weights $d_m$ can be defined such that
$$\sum_{m=0}^{r-1} d_m(x) P^{LO}_m(x) = P^{HO}(x).$$
However, the optimal linear reconstruction may suffer from Gibbs phenomena, when strong gradients are present.
To avoid this, the linear convex combination is modified by introducing the following non-linear weights:
\begin{equation}
\omega_m = \frac{\alpha_m}{\sum^{r-1}_{k=0}\alpha_k}  , \quad\text{ where } \quad 	\alpha_k = \frac{d_k}{(\beta_k+\epsilon)^2}  .
\end{equation}
In the last expression $\epsilon$ is a small number used to avoid division by zero (in general, $\epsilon=10^{-6}$),
while the $\beta_k$ are the smoothness indicators defined as
\begin{equation}
\beta_k = \sum_{j=1}^{r-1} \int_{\xin}^{\xip} \left(\frac{\diff{}^j}{\diff{x}^j} P^{LO}_k(x)\right)^2 \Delta x^{2j-1}\diff{x} , \qquad k = 0,\ldots,r-1.
\end{equation}
Finally, the WENO reconstructed polynomial is defined as:
\begin{equation}
\tilde{u}(x) = \sum_{m=0}^{r-1} \omega_m(x) P^{LO}_m(x).
\end{equation}
\begin{remark}[WENO reconstruction of global fluxes]
At steady states, the global flux $\bG$ is constant, and the high order WENO reconstruction returns the same value  at all quadrature points and at the cell interfaces.
In combination with the fact that the numerical fluxes are defined as a function of the global flux, this ensures that the numerical method is well-balanced.
\end{remark}

The values of the optimal weights $d_m$ and the formulae for computing $\beta_m$ can be found
in~\cite{jiang1996efficient,balsara2000monotonicity} up to $r=6$.
In this work, we will test the reconstruction with orders $p=1$ (standard piece-wise constant), $p=3$ and $p=5$.

\subsection{Bathymetry source term treatment} \label{sec:larSourceDiscretization}

For the standard shallow water and shallow water moment models, the bathymetry source is crucial to simulate water flows.
Exact preservation of lake at rest solutions can be directly embedded for the SWME by following the approach presented in \cite{ciallella2023arbitrary}.
In this section, we recall explicit formulas for the evaluation of the bathymetry term at quadrature points to achieve exact well-balancedness for lake at rest equilibria.
For general shallow water moment models, the lake at rest solution is simply defined as
\begin{equation}\label{eq:lar}
\eta = h + b \equiv \eta_0, \quad \text{ and } \quad u_m = \alpha_1 = \ldots = \alpha_N \equiv 0.
\end{equation}
We also discuss the issue of the jumps at the interface by focusing for simplicity on the frictionless case.

Following standard approaches to achieve exact well-balancedness at lake at rest equilibria, the reconstruction is performed on the free surface elevation $\eta$ and bathymetry $b$ (using the same WENO weights as $\eta$).
We denote the reconstructed values of $\eta$ and $b$ at the quadrature points by $\tilde{\eta}_{i,q}$ and $\tilde{b}_{i,q}$ respectively.
Then, we can define the reconstructed values of the water height $h_{i,q} = \tilde{\eta}_{i,q} -\tilde{b}_{i,q}$.
Let us also define a Lagrange interpolation of the bathymetry inside the cell $\Omega_i$ and its evaluation at the interfaces:
\begin{equation}\label{eq:bath_reconstruction}
\tilde{b}_i(x):= \sum_{q} L_q(x) \tilde{b}_{i,q},\qquad \text{and} \qquad b^L_{\iip} = \tilde{b}_i(x_{\iip}), \,\qquad b^R_{\iin} = \tilde{b}_i(x_{\iin}) .
\end{equation}
Following \cite{xing2005high}, we re-write the bathymetry term as
\begin{equation}\label{eq:sourceWB}
S_b(x) = - g h(x) \partial_x b(x) =- g \eta(x) \partial_x b(x) + g\partial_x\left(\frac{b^2(x)}{2}\right),
\end{equation}
and we compute the source integral for the momentum equation $\mathcal{R}_m$ in the quadrature point $x_{i,q}$ as
\begin{equation}\label{eq: Rquad2}
\begin{split}
	(\mathcal R_m)_{i,q} &= (\mathcal R_m)^R_{\iin} - \int_{\xin^R}^{x_{i,q}} S_b(x)\diff{x} \\
					 &= (\mathcal R_m)^R_{\iin} + g\int_{\xin^R}^{x_{i,q}} \eta(x)\partial_x b(x)  \diff{x} - g\left(\frac{(b_{i,q})^2}{2} - \frac{(b_{\iin}^R)^2}{2}\right) \\
					 &= (\mathcal R_m)^R_{\iin} + g\int_{\xin^R}^{x_{i,q}}\sum_\theta L_\theta(x) \tilde{\eta}_{i,\theta}\sum_s L'_s(x_{i,\theta}) \tilde{b}_{i,s} \diff{x} - g\left(\frac{(\tilde b_{i,q})^2}{2} - \frac{(b_{\iin}^R)^2}{2}\right),
\end{split}
\end{equation}
where in each quadrature point  $x_{i,q}$ we set
$$
\partial_x b(x_{i,q}) = \sum_s L'_s(x_{i,q}) b(x_{i,s}) .
$$
As a matter of fact, the bathymetry term and non-conservative products are treated using the same approach to reach high-order accuracy.

Thanks to \cref{eq: Rquad2}, we can define the left interface terms, introduced in \eqref{eq: Rinterface}, as
\begin{equation}\label{eq: Rinterface2}
	\begin{split}
	(\mathcal{R}_m)_{i+1/2}^L &= (\mathcal{R}_m)_{i-1/2}^R + g\int_{\xin^R}^{\xip^L} \eta(x)\partial_x b(x)  \diff{x} - g\left(\frac{(b_{\iip}^L)^2}{2} - \frac{(b_{\iin}^R)^2}{2}\right) \\
			& = (\mathcal{R}_m)_{i-1/2}^R -\Delta x(\bar{S}_b)_i,
	\end{split}
\end{equation}
where the cell average of the source term is defined as
$$ (\bar{S}_b)_i : = \frac{1}{\Delta x}\left[g\int_{\xin^R}^{\xip^L} \sum_\theta L_\theta(x) \tilde{\eta}_{i,\theta}\sum_s L'_s(x_{i,\theta}) \tilde{b}_{i,s}  \diff{x} - g\left(\frac{(b_{\iip}^L)^2}{2} - \frac{(b_{\iin}^R)^2}{2}\right)  \right]. $$
Finally, we need a recursive definition for the jump across the interfaces.
In \cite{ciallella2023arbitrary}, we proposed to define the jump of the bathymetric source term as
\begin{equation}\label{eq:jump}
[\![ \mathcal R_m ]\!]_{\iip} := g\frac{\eta_\iip^R +\eta_\iip^L }{2}\left(b_\iip^R  - b_\iip^L\right) - g\left(\frac{(b_\iip^R)^2}{2} - \frac{(b_\iip^L)^2}{2}\right).
\end{equation}
In \ref{app:lar} we recall the proof to achieve the well-balanced property of the scheme for the lake at rest equilibria,
and show that this definition of the jump is the only one that allows to achieve this property.
These  definitions allow to easily prove the following property also for general shallow water moment models.
This is due to the fact that, at lake at rest equilibria \eqref{eq:lar}, the velocity and moments are equal to zero which brings back the model
to the classical shallow water case.
\begin{proposition}[Lake at rest preservation] \label{prop:lar}
	The global flux finite volume WENO scheme with quadrature \eqref{eq: Rquad2} of the bathymetry term,
	and with the definition \eqref{eq:jump} at the interface is exactly well-balanced for the lake at rest steady state.
\begin{proof} See \ref{app:lar}.
\end{proof}
\end{proposition}
\begin{remark}[Interface jump of the bathymetry source term]\label{remark:path}
\cref{eq:jump} is similar to classical strategies used in path conservative methods \cite{CASTRO2017131},
where one uses a linear path to connect the left and right states when evaluating the integral.
However, in our case, the jump is added to the global flux that naturally entails general equilibria preservation.
As a matter of fact, the same jump can be retrieved by integrating the source term across the interface, as follows:
\begin{equation}
	\begin{split}
	[\![ \mathcal R_m ]\!]_{\iip} &= g\int_{\xip^L}^{\xip^R} \eta(x)\partial_x b(x)  \diff{x} - g\left(\frac{(b_{\iip}^R)^2}{2} - \frac{(b_{\iip}^L)^2}{2}\right) \diff{x} \\
				                &= g\int_{0}^{1} \eta(\Psi(s))\partial_s b(\Psi(s)) \diff{s} - g\left(\frac{(b_{\iip}^R)^2}{2} - \frac{(b_{\iin}^R)^2}{2}\right) \diff{x}.
	\end{split}
\end{equation}
The jump in \cref{eq:jump} can be directly computed using the linear path for $\eta$ (and similarly for $b$) defined as
$$\Psi:[0,1]\times \R \times \R \rightarrow \R, \qquad \Psi(0;\eta^L,\eta^R) = \eta^L, \qquad \Psi(0;\eta^L,\eta^R) = \eta^R .$$
For instance, for the free surface elevation $\eta$, we can define the linear path as $\eta(s) = \eta^L + s (\eta^R - \eta^L).$
This also tells us that the only path that preserves the well-balanced property is the linear one.
\end{remark}
\subsection{Non-conservative product treatment} \label{sec:nonconservative}

In this section, we will discuss the treatment of non-conservative products in the context of the global flux WENO finite volume method.
As also shown in \cref{sec:moment}, several hyperbolic systems used to describe shallow flows present non-conservative products in their formulation.
One of the goals of this work is to achieve consistency when non-conservative products are present, and corroborate the expectations
by testing the method on several complex shallow water moment models.
In particular, the main idea consists in adding the non-conservative product to the global fluxes, and treating them as an additional source term.
Then, the resulting system in the global flux formulation can be advanced in time using the numerical fluxes defined in \cref{sec_Space},
by considering the full system matrix $\bA = \partial_\bu \bF -\bB$, rather than the simple flux Jacobian.

To explain how this can be done, for the sake of simplicity, we consider the frictionless case of the first order shallow water moment (SWME1)
model described in \cref{eq:SWME1} with \eqref{eq:SWME1-B}, but clearly the same approach has been also applied to higher order models
as also shown in the numerical simulations in \cref{sec:numerics}.
Since the integral $\mathcal R_m$ can be defined as in \cref{sec:larSourceDiscretization},
we will focus on the integral $\mathcal R_1$ that concerns the non-conservative product.
In the same spirit of the bathymetry source term, we provide here the explicit formulas
of the global flux for the moment equation at the quadrature points:
\begin{equation}
	\begin{split}
		(\mathcal R_1)_{i,q} &= (\mathcal R_1)^R_{\iin} - \int_{\xin^R}^{x_{i,q}} B_{33}\partial_x U_3 \diff{x}
		 = (\mathcal R_1)^R_{\iin} - \int_{\xin^R}^{x_{i,q}} u_m \partial_x (h \alpha_1) \diff{x} \\
		&= (\mathcal R_1)^R_{\iin} - \int_{\xin^R}^{x_{i,q}}\sum_\theta L_\theta(x) (\widetilde{u_m})_{i,\theta}\sum_s L'_s(x_{i,\theta}) (\widetilde{h \alpha_1})_{i,s} \diff{x}.
	\end{split}
\end{equation}
Similarly, the left interface terms can be simply defined as
\begin{equation}
	\begin{split}
		(\mathcal R_1)^L_{\iip} &= (\mathcal R_1)^R_{\iin} - \int_{\xin^R}^{\xip^L} B_{33}\partial_x U_3 \diff{x}
		= (\mathcal R_1)^R_{\iin} - \Delta x (\overline{ B_{33}\partial_x U_3})_i,
	\end{split}
\end{equation}
where
$$ (\overline{ B_{33}\partial_x U_3})_i := \frac{1}{\Delta x}\left[\int_{\xin^R}^{\xip^L} \sum_\theta L_\theta(x) (\widetilde{u_m})_{i,\theta}\sum_s L'_s(x_{i,\theta}) (\widetilde{h \alpha_1})_{i,s} \diff{x} \right]. $$
Since a derivative is present in the non-conservative product, it becomes crucial
for consistency reasons to define a jump across the interface, such that
$$ (\mathcal R_1)^R_{\iip} = (\mathcal R_1)^L_{\iip} +  [\![ \mathcal R_1 ]\!]_{\iip} ,$$
similarly to how it would be done for path-conservative methods in standard finite volume formulation.
In this case, we can retrieve the jump by integrating the non-conservative product across the interface, as shown in \cref{remark:path}:
\begin{equation}
	\begin{split}
		[\![ \mathcal R_1 ]\!]_{\iip} &= \int_{\xip^L}^{\xip^R} u_m \partial_x(h \alpha_1) \diff{x}
				                = \int_{0}^{1} u_m(\Psi(s))\partial_s(h \alpha_1)(\Psi(s))  \diff{s} \\
								&= \frac{(u_m)^R_{\iip}+(u_m)^L_{\iip}}{2}\left((h\alpha_1)^R_{\iip} - (h\alpha_1)^L_{\iip}  \right).
	\end{split}
\end{equation}

\subsection{Additional source term treatment: friction} \label{sec:friction}

Additional source terms, such as the friction term for the general shallow water moment models,
can be easily accounted for in the global flux formulation.
For the sake of simplicity, we will focus on the friction term for the SWME1 described in \cref{eq:SWME1} with \eqref{eq:SWME1-B},
and provide the explicit formulas for the global flux in quadrature points and interfaces.
The formulas can be easily extended to higher order models, and to other source terms.
In this case, the integral terms $\bR$ in quadrature points can thus be defined as
\begin{equation}
	\begin{split}
		&(\mathcal R_m)_{i,q} = (\mathcal R_m)^R_{\iin} + \int_{\xin^R}^{x_{i,q}}\sum_\theta L_\theta(x) \left( g \tilde{\eta}_{i,\theta}\sum_s L'_s(x_{i,\theta}) \tilde{b}_{i,s} + \frac{\nu}{\lambda}\left((\widetilde{u_m})_{i,\theta} + (\widetilde{\alpha_1})_{i,\theta}\right) \right) \diff{x} \\ &\qquad\qquad -g\left(\frac{(\tilde b_{i,q})^2}{2} - \frac{(b_{\iin}^R)^2}{2}\right) \\
		&(\mathcal R_1)_{i,q} = (\mathcal R_1)^R_{\iin}  \\ &\qquad-\int_{\xin^R}^{x_{i,q}}\sum_\theta L_\theta(x) \left( (\widetilde{u_m})_{i,\theta}\sum_s L'_s(x_{i,\theta}) (\widetilde{h \alpha_1})_{i,s} + 3\frac{\nu}{\lambda} \left((\widetilde{u_m})_{i,\theta} + \left(1+4\frac{\lambda}{h}\right)(\widetilde{\alpha_1})_{i,\theta}\right) \right)\diff{x} .
	\end{split}
\end{equation}
%

\section{Time discretization}\label{sec:time}

Time integration can be performed with any time discretization method in the spirit of the method of lines, e.g., standard RK schemes.
In this work, we employ a  Deferred Correction (DeC) method because it is a family of one step methods with arbitrarily high order of accuracy.
The original DeC formulation was introduced in~\cite{daniel1968iterated}, then developed and studied in its different forms in~\cite{dutt2000dec, minion2003dec,christlieb2010integral, liu2008strong}. A
slightly different form was presented in~\cite{abgrall2017high} for applications to finite element methods.  The DeC method is presented as an iterative procedure that involves two operators.
The iteration process mimics a Picard--Lindel\"of iteration at the discrete level with a fixed-point iterative method.
Each iteration aims at gaining  one order of accuracy, so that the order of accuracy sought can be reached with a finite number of corrections.
Since the DeC method is not a part of the main contributions of this work, we will not go into details here.
The interested reader is referred to the aforementioned references for a complete description of the method, or to our previous work \cite{ciallella2023arbitrary}.
More recent developments of the DeC methods can be found in
\cite{offner2020arbitrary,veiga2021dec,abgrall2021relaxation,ciallella2021arbitrary,izgin2023study,micalizzi2024new,veiga2024improving,micalizzi2025efficient,micalizzi2024novel,ciallella2025high}.

\section{Numerical experiments}\label{sec:numerics}

In this section, we present several numerical experiments to validate the global flux WENO
finite volume method for the shallow water moment models.
In particular, the contributions of this work consider both methodology and applications:
on the one hand, we show that the new method is able to deal with general non-conservative products
in a straightforward way, and on the other hand that it is possible to build general high order
finite volume well-balanced methods for complex models, even when no closed form of the equilibria is known.
The latter feature is particularly important for systems such as the nonlinear shallow water moment models
because it allows to simulate in a much more accurate way more realistic hyperbolic systems,
without the constraint of knowing a priori the exact form of the studied problems.

For simplicity, we start by validating through convergence analysis the method for the
SWME1 model, for which an analytical steady state solution can be computed.
The exact preservation of the lake at rest solution is studied numerically, along with
classical supercritical and subcritical moving equilibria.
A convergence analysis is shown on a set of uniform meshes of $N_e=100,200,400,600,800$ control volumes,
with WENO1 (standard piecewise constant), WENO3 and WENO5 reconstructions.
Notice that the name WENO1 to define standard piecewise constant reconstructions is given only to make it easier to understand method acronyms.
Since for the SWME1 model it is possible to easily find eigenvalues and eigenvectors,
we performed the convergence analysis with both numerical global fluxes in \cref{eq:upwind flux} and \cref{eq:central flux}.
When dealing with more complex models, we will focus simply on the one in \cref{eq:central flux}, which does not require any decomposition.
The latter simplifies not only the implementation for complex models, due to much easier formulation,
but it also decreases notably the computational cost of the method, which is crucial when increasing the complexity of the studied models.
Moreover, rather than focusing on the comparisons between the new well-balanced methods and non-well-balanced ones,
which has already been thoroughly studied in \cite{ciallella2023arbitrary} for the classical SWE,
the rest of the numerical experiments will focus on comparing the numerical results obtained with several
shallow water moment models.
The goal is to show that the method achieves comparable accuracy with both simplified and fully nonlinear moment models,
and therefore the new well-balanced methods can be used to study much more complex models never tackled before with
other well-balanced methods.
For the simulation of steady states with friction, we take $\nu=0.05$ and $\lambda=1$, similarly to what was used in \cite{pimentel2022fully}.

\subsection{Lake at rest and its perturbation}

The first test case is the lake at rest solution characterized by the following initial and exact conditions for the SWME1 model:
$$ h(x,0) = \eta_0 - b(x),\qquad u_m = \alpha_1 = 0, \quad\text{ with }\quad \eta_0 \equiv 1 $$
in a computational domain of size $[0,25]$ with subcritical inlet/outlet conditions at the left and right boundaries.
The bathymetry considered is $\mathcal C^\infty$, with values smaller than machine precision at boundaries, and reads
\begin{equation}\label{eq:bathy}
	b(x) = 0.05 \sin(x-x_0)\exp(1-(x-x_0)^2),\quad\text{ with }\quad x_0 =12.5.
\end{equation}
The gravitational constant is set to $g=1$ and all simulations are run until final time $T=1$.
To make the equation on the moment $h \alpha_1$ non-trivial, we consider for all these simulations
the friction source term.
Considering only the bathymetry term would yield $\alpha_1 \equiv 0$, when starting from a zero initial condition (which is the case for the lake at rest).

In \cref{tab:LAR_UPWIND} and \cref{tab:LAR_SUPG}, we present the discretization errors obtained with the three
WENO1, WENO3 and WENO5 reconstructions, coupled with the numerical global fluxes $\widehat \bG^u$ and $\widehat \bG^c$, presented in \cref{eq:upwind flux}
and \cref{eq:central flux} respectively. It can be noticed that for all reconstructions, and all levels of refinement,
the methods preserve the solution up to machine precision, confirming the theoretical expectations.

\begin{table}
	\caption{Lake at rest for SWME1: errors and estimated order of accuracy (EOA) with the global flux schemes, using WENO1, WENO3 and WENO5 reconstructions and numerical flux $\widehat \bG^u$.}\label{tab:LAR_UPWIND}
	\scriptsize
	\centering
	\begin{tabular}{c|cc|cc|cc} \hline\hline
			&\multicolumn{2}{c|}{$h$} &\multicolumn{2}{c|}{$hu_m$} &\multicolumn{2}{c}{$h\alpha_1$}  \\[0.5mm]
			\cline{2-7}
			$N_e$ & $L_2$ error        & EOA & $L_2$  error       & EOA & $L_2$ error       &EOA  \\ \hline \hline
			&\multicolumn{6}{c}{GF-WENO1}\\ \hline
			100 &    0.000   & -- &     1.324E-017  & -- &  8.295E-019 & -- \\
			200 &    0.000   & -- &     1.840E-017  & -- &  1.128E-018 & -- \\
			400 &    0.000   & -- &     1.607E-017  & -- &  1.096E-018 & -- \\
			600 &    0.000   & -- &     2.427E-017  & -- &  1.826E-018 & -- \\
			800 &    0.000   & -- &     2.333E-017  & -- &  1.720E-018 & -- \\ \hline    		
			&\multicolumn{6}{c}{GF-WENO3}\\ \hline
			100 &    0.000      & -- &  1.074E-016 & -- &  5.900E-018 & -- \\
			200 &    0.000      & -- &  1.262E-016 & -- &  7.664E-018 & -- \\
			400 &    5.551E-019 & -- &  1.723E-016 & -- &  1.151E-017 & -- \\
			600 &    4.625E-018 & -- &  2.012E-016 & -- &  1.444E-017 & -- \\
			800 &    4.718E-018 & -- &  2.458E-016 & -- &  1.750E-017 & -- \\   \hline
			&\multicolumn{6}{c}{GF-WENO5}\\  \hline
			100 &    2.220E-018 & -- &  1.072E-016 & -- &  6.779E-018 & -- \\
			200 &    4.996E-018 & -- &  1.128E-016 & -- &  7.421E-018 & -- \\
			400 &    1.332E-017 & -- &  1.971E-016 & -- &  1.400E-017 & -- \\
			600 &    2.646E-017 & -- &  2.281E-016 & -- &  1.582E-017 & -- \\
			800 &    2.386E-017 & -- &  2.288E-016 & -- &  1.664E-017 & -- \\   \hline  			
			\hline
	\end{tabular}
\end{table}

\begin{table}
	\caption{Lake at rest for SWME1: errors and estimated order of accuracy (EOA) with the global flux schemes, using WENO1, WENO3 and WENO5 reconstructions and numerical flux $\widehat \bG^c$}\label{tab:LAR_SUPG}
	\scriptsize
	\centering
	\begin{tabular}{c|cc|cc|cc} \hline\hline
			&\multicolumn{2}{c|}{$h$} &\multicolumn{2}{c|}{$hu_m$} &\multicolumn{2}{c}{$h\alpha_1$}  \\[0.5mm]
			\cline{2-7}
			$N_e$ & $L_2$ error        & EOA & $L_2$  error       & EOA & $L_2$ error       &EOA  \\ \hline \hline
			&\multicolumn{6}{c}{GF-WENO1}\\ \hline
			100 &    0.000   & -- &     3.028E-018 & -- &  1.999E-019 & -- \\
			200 &    0.000   & -- &     8.607E-018 & -- &  5.914E-019 & -- \\
			400 &    0.000   & -- &     4.276E-018 & -- &  3.283E-019 & -- \\
			600 &    0.000   & -- &     1.066E-017 & -- &  8.027E-019 & -- \\
			800 &    0.000   & -- &     7.421E-018 & -- &  6.282E-019 & -- \\    \hline 		
			&\multicolumn{6}{c}{GF-WENO3}\\ \hline
			100 &    0.000      & -- &  7.247E-017 & --&  4.587E-018 & -- \\
   			200 &    0.000      & -- &  1.035E-016 & --&  6.699E-018 & -- \\
   			400 &    0.000      & -- &  1.239E-016 & --&  8.196E-018 & -- \\
   			600 &    9.251E-019 & -- &  1.740E-016 & --&  1.220E-017 & -- \\
   			800 &    8.326E-019 & -- &  2.274E-016 & --&  1.552E-017 & -- \\    \hline
			&\multicolumn{6}{c}{GF-WENO5}\\  \hline
			100 &    1.110E-018 & -- &  9.246E-017 & -- &  5.960E-018 & -- \\
			200 &    2.775E-018 & -- &  1.146E-016 & -- &  7.393E-018 & -- \\
			400 &    1.249E-017 & -- &  1.895E-016 & -- &  1.304E-017 & -- \\
			600 &    2.553E-017 & -- &  2.092E-016 & -- &  1.547E-017 & -- \\
			800 &    2.498E-017 & -- &  2.081E-016 & -- &  1.478E-017 & -- \\    \hline 		
			\hline
	\end{tabular}
\end{table}

\begin{figure}
	\centering
	\subfigure{
		\includegraphics{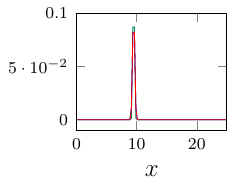}}\quad
	\subfigure{
		\includegraphics{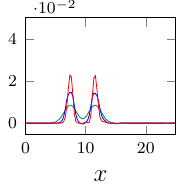}}\quad
	\subfigure{
		\includegraphics{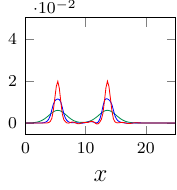}}\quad
	\subfigure{
		\includegraphics{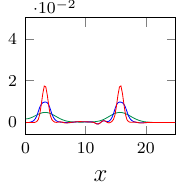}}\quad
	\subfigure{
		\includegraphics{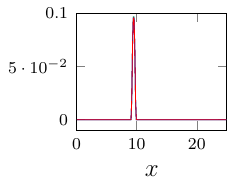}}\quad
	\subfigure{
		\includegraphics{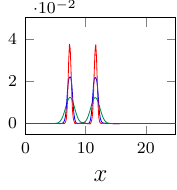}}\quad
	\subfigure{
		\includegraphics{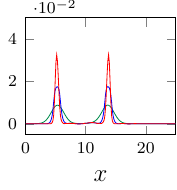}}\quad
	\subfigure{
		\includegraphics{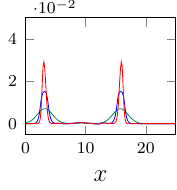}}\quad
				\setcounter{subfigure}{0}
	\subfigure[$t=0$]{
		\includegraphics{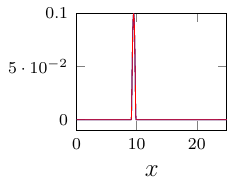}}\quad
	\subfigure[$t=0.66$]{
		\includegraphics{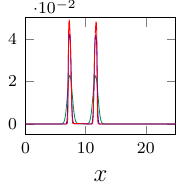}}\quad
	\subfigure[$t=1.33$]{
		\includegraphics{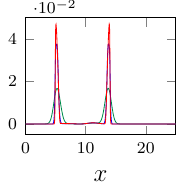}}\quad
	\subfigure[$t=2$]{
		\includegraphics{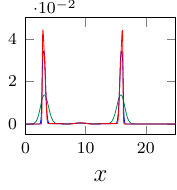}}
	\caption{Small perturbation of lake at rest computed with GF-WENO1 (green), GF-WENO3 (blue) and GF-WENO5 (red): value of $h-h_{eq}$ with $N_e=100$ (top), $N_e=200$ (middle), $N_e=800$ (bottom) at different simulation times.}\label{fig:pert_lar_h}
	\end{figure}

\begin{figure}
		\centering
	\subfigure{
		\includegraphics{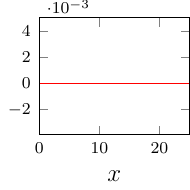}}\quad
	\subfigure{
		\includegraphics{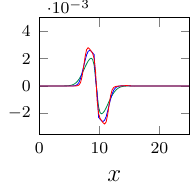}}\quad
	\subfigure{
		\includegraphics{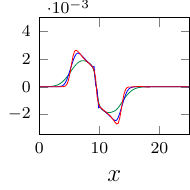}}\quad
	\subfigure{
		\includegraphics{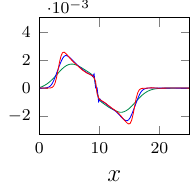}}\quad
	\subfigure{
		\includegraphics{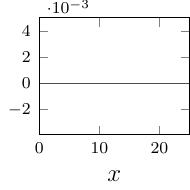}}\quad
	\subfigure{
		\includegraphics{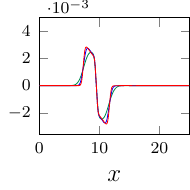}}\quad
	\subfigure{
		\includegraphics{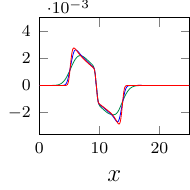}}\quad
	\subfigure{
		\includegraphics{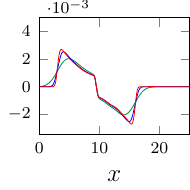}}\quad
	\setcounter{subfigure}{0}
	\subfigure[$t=0$]{
		\includegraphics{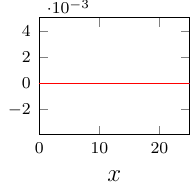}}\quad
	\subfigure[$t=0.66$]{
		\includegraphics{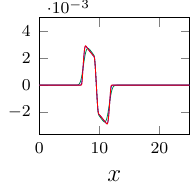}}\quad
	\subfigure[$t=1.33$]{
		\includegraphics{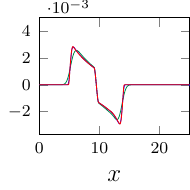}}\quad
	\subfigure[$t=2$]{
		\includegraphics{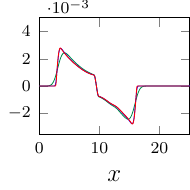}}
	\caption{Small perturbation of lake at rest computed with GF-WENO1 (green), GF-WENO3 (blue) and GF-WENO5 (red): value of $(h\alpha_1)-(h\alpha_1)_{eq}$ with $N_e=100$ (top), $N_e=200$ (middle), $N_e=800$ (bottom) at different simulation times.}\label{fig:pert_lar_ha1}
	\end{figure}

In \cref{fig:pert_lar_h} and \cref{fig:pert_lar_ha1}, we present a perturbation analysis for the lake at rest steady
state with a perturbation on the water height of the following shape:
$$\delta h(x,0) = 10^{-1} \exp\left(1 - \frac{1}{(1-r(x))^2}\right),\quad\text{ with }\quad r(x) = 4(x-9.5)^2. $$
The simulation is run starting from the equilibria solution $h_{eq}$, by adding the perturbation to the water height as
$h(x,0) = h_{eq}(x,0) + \delta h(x,0)$.

We present the numerical results obtained with $g=9.8$ on different levels of refinement to show the qualitative convergence of the solution.
Larger gravity values are used here simply to have a more pronounced perturbation effect.
\cref{fig:pert_lar_ha1} also shows a variation of $\alpha_1$ to due the friction term. This small variation of the order of $\sim 10^{-3}$ can be studied accurately
with the new methods. Also in \cref{fig:pert_lar_h} we can see the effect of friction on the water height
perturbation, which decreases with time.

\subsection{Steady states without friction and their perturbation}

In this section, we test the numerical methods on the SWME1 model for classical steady states without friction.
We stick to the SWME1 model because thanks to its simple formulation, we can compute analytically the exact moving
equilibria to perform convergence analysis of the solution.
For the same reason, we consider the same bathymetry defined in \cref{eq:bathy} which is smooth enough to be tested in
the context of high order schemes.
We consider the following two sets of final time $T$ (to reach steady state, i.e.\ time residual up to machine precision),
initial and boundary conditions:
\begin{itemize}
	\item Supercritical case
		\begin{alignat*}{3}
			&T = 50, \\
			&h(x,0) = 2 - b(x),\qquad &&h u_m (x,0) = 0  ,\qquad && h \alpha_1 (x,0) = -0.5, \\
			&h(0,t) = 2       ,\qquad &&h u_m (0,t) = 24 ,\qquad && h \alpha_1 (0,t) = -0.5. \\
		\end{alignat*}
	\item Subcritical case
	\begin{alignat*}{3}
		&T = 400, \\
		&h(x,0) = 2 - b(x),\qquad &&h u_m (x,0) = 0  ,\qquad && h \alpha_1 (x,0) = 0.1, \\
		&h(25,t) = 2       ,\qquad &&h u_m (0,t) = 4.42 ,\qquad && h \alpha_1 (0,t) = 0.1. \\
	\end{alignat*}
\end{itemize}

Following also other references \cite{cheng2019new}, for these tests, the gravitational constant is set to $g=9.812$,
and the convergence analysis is performed by comparing the numerical solution to the analytical one
computed by solving \cref{eq:SWME1-exactH} with a nonlinear solver.

\begin{figure}
	\centering
	\subfigure[Free surface $\eta$]{
		\includegraphics{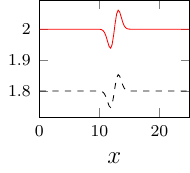}}\quad
	\subfigure[Average speed $u_m$]{
		\includegraphics{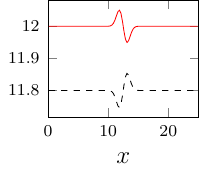}}\quad
	\subfigure[Moment $\alpha_1$]{
		\includegraphics{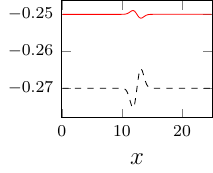}}\quad
	\subfigure[Distribution $u(\zeta)$]{
		\includegraphics{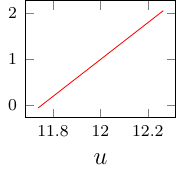}}
	\caption{Supercritical case without friction computed with GF-WENO5 for the SWME1 model: characteristic variables (red) and rescaled bathymetry (black dashed) with $N_e=100$.}\label{fig:super_var}
\end{figure}

\begin{figure}
	\centering
	\subfigure[Free surface $\eta$]{
		\includegraphics{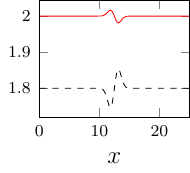}}\quad
	\subfigure[Average speed $u_m$]{
		\includegraphics{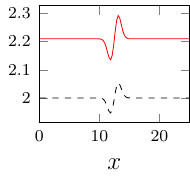}}\quad
	\subfigure[Moment $\alpha_1$]{
		\includegraphics{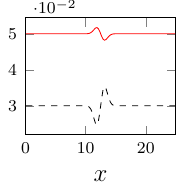}}\quad
	\subfigure[Distribution $u(\zeta)$]{
		\includegraphics{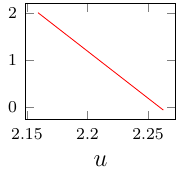}}	
	\caption{Subcritical case without friction computed with GF-WENO5 for the SWME1 model: characteristic variables (red) and rescaled bathymetry (black dashed)  with $N_e=100$.}\label{fig:sub_var}
\end{figure}

In \cref{fig:super_var} and \cref{fig:sub_var}, we present the characteristic variables of both supercritical and subcritical equilibria
obtained with GF-WENO5 on a mesh of $N_e=100$ elements.
In particular, we show the behaviors obtained at steady state for the free surface elevation $\eta$, the average speed $u_m$,
the first order moment $\alpha_1$ and the velocity distribution along the vertical direction computed at randomly chosen point $x=18.62$ with \cref{eq:veldistribution} written as
$$	u(\zeta) = u_m + \left(1-2\zeta\right)\alpha_1 $$
where $\zeta=\frac{z-b}{h}$.
In the supercritical case, we observe that the free surface elevation $\eta$ follows the bathymetry $b$,
while the average speed $u_m$ grows as $b$ decreases and descreses as $b$ increases.
As expected, the opposite is observed for the subcritical case.
Due to the choice of the first order moment $\alpha_1$, the velocity distribution along the vertical direction is linear,
with a slope that depends on the value of $\alpha_1$.

In \cref{tab:SUPER_UPWIND} and \cref{tab:SUPER_SUPG}, we show the convergence analysis performed on the supercritical case
with all reconstructions and the numerical global fluxes presented in \cref{eq:upwind flux} and \cref{eq:central flux}, respectively.
While in \cref{tab:SUB_UPWIND} and \cref{tab:SUB_SUPG}, we present the same results obtained for the subcritical case.
Estimated order of accuracy for all cases highlights that the new method is able to achieve the expected high order convergence properties.
In particular, we would like to point out that although the piecewise constant reconstruction (GF-WENO1) is in general
first order accurate, a superconvergence to second order is achieved for stationary problems. However, when non-stationary
cases are considered, e.g., perturbed steady states, the method is indeed only first order accurate, which can be already noticed
from \cref{fig:pert_lar_h}.

\begin{table}
	\caption{Supercritical steady state for SWME1: errors and estimated order of accuracy (EOA) with the global flux schemes, using WENO1, WENO3 and WENO5 reconstructions and numerical flux $\widehat \bG^u$.}\label{tab:SUPER_UPWIND}
	\scriptsize
	\centering
	\begin{tabular}{c|cc|cc|cc} \hline\hline
			&\multicolumn{2}{c|}{$h$} &\multicolumn{2}{c|}{$hu_m$} &\multicolumn{2}{c}{$h\alpha_1$}  \\[0.5mm]
			\cline{2-7}
			$N_e$ & $L_2$ error        & EOA & $L_2$  error       & EOA & $L_2$ error       &EOA  \\ \hline \hline
			&\multicolumn{6}{c}{GF-WENO1}\\ \hline
			100  & 8.424e-06  &  --   & 2.096e-14 & --  & 4.214e-06 &  --   \\
			200  & 2.133e-06  & 1.98  & 1.698e-14 & --  & 1.067e-06 & 1.98  \\
			400  & 5.321e-07  & 2.00  & 5.854e-14 & --  & 2.662e-07 & 2.00  \\
			600  & 2.364e-07  & 2.00  & 1.987e-14 & --  & 1.182e-07 & 2.00  \\
			800  & 1.329e-07  & 2.00  & 2.305e-14 & --  & 6.652e-08 & 2.00  \\ \hline
			&\multicolumn{6}{c}{GF-WENO3}\\ \hline
			100  & 5.616e-08  &  --   & 1.268e-14 & --  & 2.401e-05 &  --  \\
			200  & 1.660e-08  & 1.76  & 1.061e-13 & --  & 1.172e-05 & 1.03 \\
			400  & 2.058e-09  & 3.01  & 5.255e-14 & --  & 1.485e-06 & 2.98 \\
			600  & 5.990e-10  & 3.04  & 4.780e-14 & --  & 4.338e-07 & 3.03 \\
			800  & 2.680e-10  & 2.80  & 1.610e-13 & --  & 1.960e-07 & 2.76 \\ \hline 		
			&\multicolumn{6}{c}{GF-WENO5}\\  \hline
			100  & 8.482e-09  &  --   & 2.259e-14 & --  & 3.751e-06 &  --  \\
			200  & 2.666e-10  & 4.99  & 3.595e-14 & --  & 1.492e-07 & 4.65 \\
			400  & 1.007e-11  & 4.73  & 4.392e-14 & --  & 6.746e-09 & 4.47 \\
			600  & 1.100e-12  & 5.46  & 6.268e-14 & --  & 7.552e-10 & 5.40 \\
			800  & 2.147e-13  & 5.68  & 7.541e-14 & --  & 1.479e-10 & 5.67 \\\hline
			\hline
	\end{tabular}
\end{table}

\begin{table}
	\caption{Supercritical steady state for SWME1: errors and estimated order of accuracy (EOA) with the global flux schemes, using WENO1, WENO3 and WENO5 reconstructions and numerical flux $\widehat \bG^c$.}\label{tab:SUPER_SUPG}
	\scriptsize
	\centering
	\begin{tabular}{c|cc|cc|cc} \hline\hline
			&\multicolumn{2}{c|}{$h$} &\multicolumn{2}{c|}{$hu_m$} &\multicolumn{2}{c}{$h\alpha_1$}  \\[0.5mm]
			\cline{2-7}
			$N_e$ & $L_2$ error        & EOA & $L_2$  error       & EOA & $L_2$ error       &EOA  \\ \hline \hline
			&\multicolumn{6}{c}{GF-WENO1}\\ \hline
			100  &  8.424e-06  &  --   & 5.684e-15 & --  & 4.214e-06  &  --  \\
			200  &  2.133e-06  &  1.98 & 2.637e-14 & --  & 1.067e-06  & 1.98 \\
			400  &  5.321e-07  &  2.00 & 7.958e-15 & --  & 2.662e-07  & 2.00 \\
			600  &  2.364e-07  &  2.00 & 3.183e-14 & --  & 1.182e-07  & 2.00 \\
			800  &  1.329e-07  &  2.00 & 1.203e-14 & --  & 6.652e-08  & 2.00 \\ \hline
			&\multicolumn{6}{c}{GF-WENO3}\\ \hline
			100  &  5.620e-08  &  --   & 5.832e-10 & --  & 2.400e-05  &  --   \\
			200  &  1.662e-08  & 1.76  & 2.938e-10 & --  & 1.171e-05  & 1.03  \\
			400  &  2.059e-09  & 3.01  & 2.093e-11 & --  & 1.484e-06  & 2.98  \\
			600  &  5.992e-10  & 3.05  & 3.179e-12 & --  & 4.338e-07  & 3.03  \\
			800  &  2.681e-10  & 2.80  & 9.050e-13 & --  & 1.960e-07  & 2.76  \\ \hline
			&\multicolumn{6}{c}{GF-WENO5}\\  \hline
			100  &  8.479e-09  &  --   & 8.783e-11 & --  & 3.752e-06  &  --  \\
			200  &  2.665e-10  & 4.99  & 1.859e-12 & --  & 1.493e-07  & 4.65 \\
			400  &  1.007e-11  & 4.73  & 9.206e-14 & --  & 6.747e-09  & 4.47 \\
			600  &  1.102e-12  & 5.46  & 5.370e-14 & --  & 7.553e-10  & 5.40 \\
			800  &  2.061e-13  & 5.83  & 1.811e-13 & --  & 1.479e-10  & 5.67 \\ \hline
			\hline
	\end{tabular}
\end{table}

\begin{table}
	\caption{Subcritical steady state for SWME1: errors and estimated order of accuracy (EOA) with the global flux schemes, using WENO1, WENO3 and WENO5 reconstructions and numerical flux $\widehat \bG^u$.}\label{tab:SUB_UPWIND}
	\scriptsize
	\centering
	\begin{tabular}{c|cc|cc|cc} \hline\hline
			&\multicolumn{2}{c|}{$h$} &\multicolumn{2}{c|}{$hu_m$} &\multicolumn{2}{c}{$h\alpha_1$}  \\[0.5mm]
			\cline{2-7}
			$N_e$ & $L_2$ error        & EOA & $L_2$  error       & EOA & $L_2$ error       &EOA  \\ \hline \hline
			&\multicolumn{6}{c|}{GF-WENO1}\\ \hline
			100 & 7.195e-05 &  --  & 4.508e-14 & -- & 3.648e-05  &  --   \\
			200 & 1.825e-05 & 1.98 & 2.789e-14 & -- & 9.242e-06  & 1.98  \\
			400 & 4.552e-06 & 2.00 & 1.407e-13 & -- & 2.305e-06  & 2.00 \\
			600 & 2.022e-06 & 2.00 & 3.715e-13 & -- & 1.024e-06  & 2.00 \\
			800 & 1.137e-06 & 2.00 & 2.097e-13 & -- & 5.761e-07  & 2.00  \\\hline
			&\multicolumn{6}{c|}{GF-WENO3}\\ \hline
			100 & 9.071e-07 &  --  & 9.201e-10 & -- & 2.961e-05 &  --  \\
			200 & 1.208e-07 & 2.91 & 9.760e-11 & -- & 2.422e-05 & 0.29 \\
			400 & 1.883e-08 & 2.68 & 1.140e-11 & -- & 4.580e-06 & 2.40 \\
			600 & 4.358e-09 & 3.61 & 2.532e-12 & -- & 1.277e-06 & 3.15 \\
			800 & 1.527e-09 & 3.65 & 2.042e-12 & -- & 5.110e-07 & 3.18 \\\hline
			&\multicolumn{6}{c|}{GF-WENO5}\\  \hline
			100 & 1.229e-07 &  --  & 5.885e-11 & -- & 9.972e-06 &  --  \\
			200 & 2.583e-09 & 5.57 & 6.563e-13 & -- & 1.856e-07 & 5.75 \\
			400 & 4.001e-11 & 6.01 & 6.342e-14 & -- & 1.728e-08 & 3.42 \\
			600 & 5.703e-12 & 4.80 & 3.882e-14 & -- & 2.576e-09 & 4.69 \\
			800 & 1.251e-12 & 5.27 & 4.191e-14 & -- & 5.697e-10 & 5.24 \\\hline		
			\hline
	\end{tabular}
\end{table}

\begin{table}
	\caption{Subcritical steady state for SWME1: errors and estimated order of accuracy (EOA) with the global flux schemes, using WENO1, WENO3 and WENO5 reconstructions and numerical flux $\widehat \bG^c$.}\label{tab:SUB_SUPG}
	\scriptsize
	\centering
	\begin{tabular}{c|cc|cc|cc} \hline\hline
			&\multicolumn{2}{c|}{$h$} &\multicolumn{2}{c|}{$hu_m$} &\multicolumn{2}{c}{$h\alpha_1$}  \\[0.5mm]
			\cline{2-7}
			$N_e$ & $L_2$ error        & EOA & $L_2$  error       & EOA & $L_2$ error       &EOA  \\ \hline \hline
			&\multicolumn{6}{c|}{GF-WENO1}\\ \hline
			100 & 7.223e-05 &  --  & 7.046e-14 & -- & 7.331e-06 &  --  \\
			200 & 1.832e-05 & 1.98 & 1.570e-13 & -- & 1.856e-06 & 1.98 \\
			400 & 4.570e-06 & 2.00 & 1.636e-13 & -- & 4.632e-07 & 2.00 \\
			600 & 2.030e-06 & 2.00 & 2.938e-13 & -- & 2.058e-07 & 2.00 \\
			800 & 1.142e-06 & 2.00 & 3.814e-13 & -- & 1.157e-07 & 2.00 \\ \hline					
			&\multicolumn{6}{c|}{GF-WENO3}\\ \hline
			100 & 8.989e-07 &  --  & 8.691e-10 & -- & 2.962e-05 &  --   \\
			200 & 1.242e-07 & 2.86 & 9.520e-11 & -- & 2.423e-05 & 0.29  \\
			400 & 1.937e-08 & 2.68 & 8.671e-12 & -- & 4.580e-06 & 2.40  \\
			600 & 4.500e-09 & 3.60 & 2.394e-12 & -- & 1.277e-06 & 3.15  \\
			800 & 1.577e-09 & 3.64 & 9.772e-13 & -- & 5.111e-07 & 3.18  \\ \hline		
			&\multicolumn{6}{c|}{GF-WENO5}\\  \hline
			100 & 1.228e-07 &  --  & 9.483e-11 & -- &  9.971e-06 &  --  \\
			200 & 2.582e-09 & 5.57 & 9.109e-13 & -- &  1.856e-07 & 5.75 \\
			400 & 3.996e-11 & 6.01 & 1.637e-13 & -- &  1.728e-08 & 3.42 \\
			600 & 5.711e-12 & 4.80 & 1.910e-13 & -- &  2.576e-09 & 4.69 \\
			800 & 1.262e-12 & 5.25 & 9.067e-14 & -- &  5.697e-10 & 5.24 \\ \hline		
			\hline
	\end{tabular}
\end{table}

\begin{figure}
	\centering
	\subfigure[Supercritical case]{
		\includegraphics{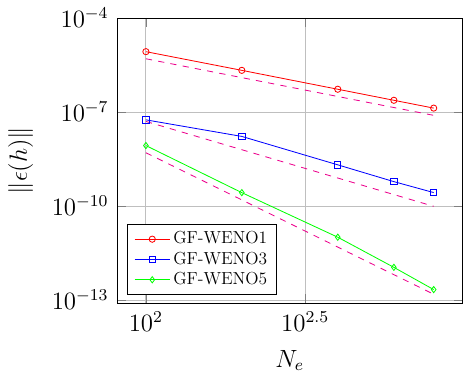}}
	\subfigure[Subcritical case]{
		\includegraphics{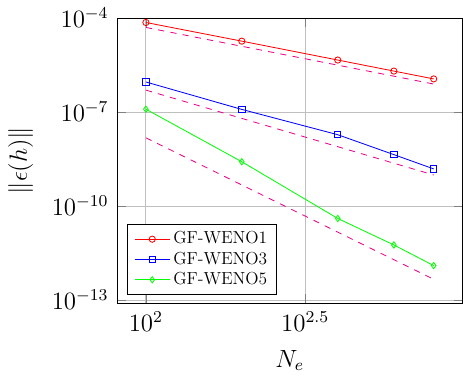}}
	\caption{Convergence analysis for supercritical (left) and subcritical (right) moving equilibria obtained for the SWME1 model. Numerical results are presented for GF-WENO1, GF-WENO3 and GF-WENO5 methods and compared against second-order, third-order, fifth-order reference curves (dashed magenta).}
	\label{fig:convplot}
	\end{figure}

\begin{figure}
		\centering
	\subfigure{
		\includegraphics{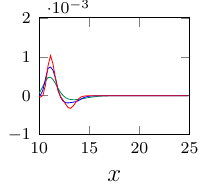}}\quad
	\subfigure{
		\includegraphics{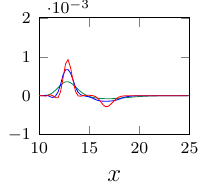}}\quad
	\subfigure{
		\includegraphics{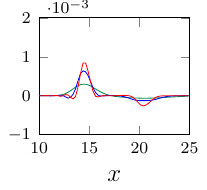}}\quad
	\subfigure{
		\includegraphics{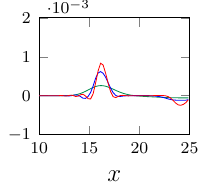}}\quad
	\subfigure{
		\includegraphics{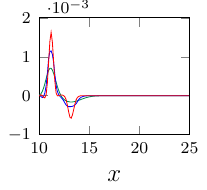}}\quad
	\subfigure{
		\includegraphics{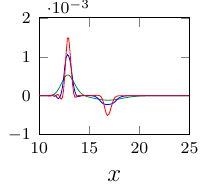}}\quad
	\subfigure{
		\includegraphics{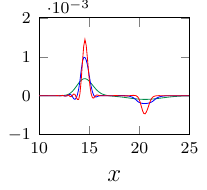}}\quad
	\subfigure{
		\includegraphics{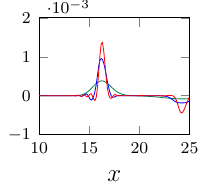}}\quad
	\setcounter{subfigure}{0}
	\subfigure[$t=0.225$]{
		\includegraphics{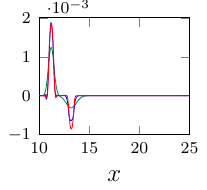}}\quad
	\subfigure[$t=0.45$]{
		\includegraphics{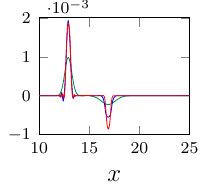}}\quad
	\subfigure[$t=0.675$]{
		\includegraphics{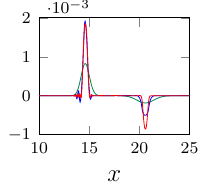}}\quad
	\subfigure[$t=0.9$]{
		\includegraphics{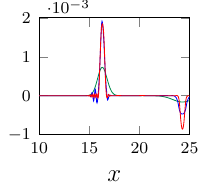}}\quad
	\caption{Small perturbation of supercritical flow without friction computed with GF-WENO1 (green), GF-WENO3 (blue) and GF-WENO5 (red): value of $h-h_{eq}$ with $N_e=100$ (top), $N_e=200$ (middle), $N_e=800$ (bottom) at different simulation times.}\label{fig:pert_super_h}
	\end{figure}

Data from the convergence tables are plotted in \cref{fig:convplot}, where discretization errors for the water height $h$ are presented on a $\log$-$\log$
scale to stress the advantages of using higher order methods.
It can be also noticed that the discretization errors obtained for the moment conservative variable, on the chosen uniform meshes, are particularly low
for the GF-WENO1 method, which achieves results comparable with GF-WENO3 even if the latter is third order accurate.
However, for either finer meshes or non-stationary problems, where the GF-WENO1 convergence order drops to one, we expect much improved results for the GF-WENO3 method,
as can be corroborated from \cref{fig:pert_lar_h} and \cref{fig:pert_super_h}.
In any case, the GF-WENO5 is by far the best and least dissipative method achieving extremely low errors for all variables.
In particular, it can be noticed that the discretization errors obtained with the new global flux method for the SWME1 model are of the same order of magnitude
as those obtained for the standard SWE model \cite{ciallella2023arbitrary}, which were several orders of magnitude lower that those obtained with standard WENO methods.
Another feature of these numerical experiments concerns the proper convergence to steady states of the simulations,
and the discretization errors computed for the average momentum $h u_m$, which should be close to machine precision.
In particular, all simulations are well converged to steady state at final time, and $h u_m$ has discretization errors
always very close to machine precision.
However, depending on the mesh refinement level, reconstruction order, numerical flux and time integration solver,
we might have that this value is $\sim 10^{-14}$ or $\sim 10^{-11}$. We would like to stress that such a low level
of errors does not influence at all the numerical results, but it emphasizes that all these parameters
may play a role in the proper convergence of the solution to the discrete steady state.

As already mentioned above, herein we only focus on the global flux methods applied to several shallow water moment models.
For detailed comparisons between the present GF-WENO methods and standard WENO methods the reader is refer to \cite{ciallella2023arbitrary}.

We also present a perturbation analysis of the supercritical case for SWME1 obtained by adding on top of the discrete steady state equilibrium a perturbation on the water height of the following shape:
$$\delta h(x,0) = 10^{-3} \exp\left(1 - \frac{1}{(1-r(x))^2}\right),\quad\text{ with }\quad r(x) = 4(x-9.5)^2. $$
As done above, the simulation is run starting from the equilibria solution $h_{eq}$, by adding the perturbation to the water height as
$h(x,0) = h_{eq}(x,0) + \delta h(x,0)$.

\begin{figure}
			\centering
	\subfigure{
		\includegraphics{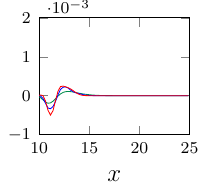}}\quad
	\subfigure{
		\includegraphics{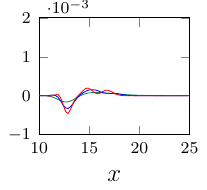}}\quad
	\subfigure{
		\includegraphics{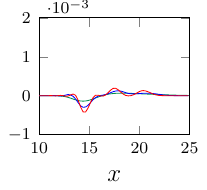}}\quad
	\subfigure{
		\includegraphics{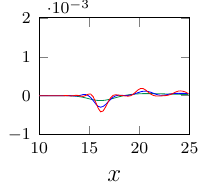}}\quad
	\subfigure{
		\includegraphics{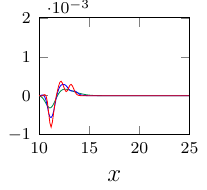}}\quad
	\subfigure{
		\includegraphics{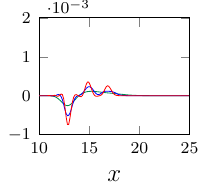}}\quad
	\subfigure{
		\includegraphics{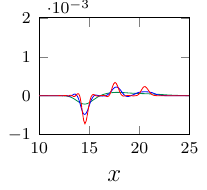}}\quad
	\subfigure{
		\includegraphics{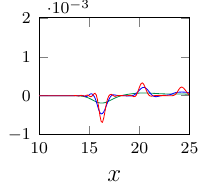}}\quad
	\setcounter{subfigure}{0}
	\subfigure[$t=0.225$]{
		\includegraphics{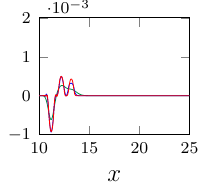}}\quad
	\subfigure[$t=0.45$]{
		\includegraphics{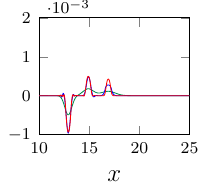}}\quad
	\subfigure[$t=0.675$]{
		\includegraphics{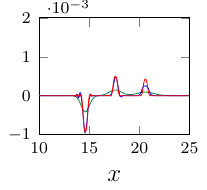}}\quad
	\subfigure[$t=0.9$]{
		\includegraphics{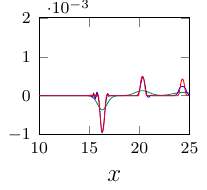}}\quad
	\caption{Small perturbation of supercritical flow without friction computed with GF-WENO1 (green), GF-WENO3 (blue) and GF-WENO5 (red): value of $(h\alpha_1)-(h\alpha_1)_{eq}$ with $N_e=100$ (top), $N_e=200$ (middle), $N_e=800$ (bottom) at different simulation times.}\label{fig:pert_super_ha1}
	\end{figure}

Since similar results can be obtained for the perturbation of the subcritical case, we do not include it for brevity.
In \cref{fig:pert_super_h} and \cref{fig:pert_super_ha1}, we show the numerical solution (as deviation from the the equilibrium solution) obtained at different times with the WENO1, WENO3, and WENO5 reconstructions
on different levels of mesh refinement. As discussed above, for both variables $h$ and $h\alpha_1$, the order of the reconstruction has a huge impact on the prediction of the perturbation evolution.
The use of three different levels of refinement also allows us to show best the high order convergence of the solution.
As expected, the GF-WENO5 method achieves the best results sharply capturing the waves appearing with the perturbation.

\subsection{Steady states with friction and their perturbation}

In the previous section, the goal was to assess the convergence properties of the new global flux method on the simplest shallow water moment model, i.e., SWME1,
due to the possibility of computing exactly a reference steady state.
In this section, we highlight the capability and generality of our approach by testing and comparing representatively four different shallow water moment models:
\begin{itemize}
\item first order shallow water moment (SWME1);
\item second order shallow water linearized moment (SWLME2);
\item second order hyperbolic shallow water moment (HSWME2);
\item second order shallow water moment (SWME2).
\end{itemize}
The challenge of this part concerns the fact that, for some of the models, it is not possible to compute equilibrium variables or exact steady states
which makes them harder to use in the context of existing well-balanced methods \cite{michel2016well,koellermeier2022steady,cao2025flux}. In particular, the generality of our global flux formulation allows us to circumvent this problem and
deal with such models no different than with the others. Moreover, notice that although herein we focus on the general shallow water moment models,
the present method is general enough to be applied to other kinds of hyperbolic balance laws with non-conservative products.

\begin{figure}
	\centering
	\subfigure[Free surface $\eta$]{
	\includegraphics{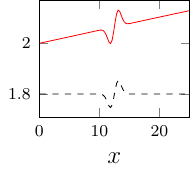}}\quad
	\subfigure[Average speed $u_m$]{
	\includegraphics{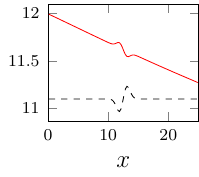}}\quad
	\subfigure[Moment $\alpha_1$]{
	\includegraphics{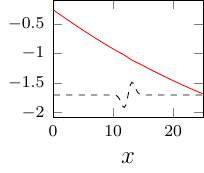}}\quad
	\subfigure[Distribution $u(\zeta)$]{
	\includegraphics{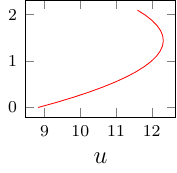}}
	\caption{Supercritical case with friction computed with GF-WENO5 for the SWME2 model: characteristic variables (red) and rescaled bathymetry (black dashed) with $N_e=100$.}\label{fig:super_SWME2}
\end{figure}

For this case with $N=2$, we include both bathymetry and friction source terms and analyze the perturbation evolution for the different models.
\cref{fig:super_SWME2} shows the characteristic variables obtained for a supercritical steady state with the GF-WENO5 method, which will be the starting point
for the perturbation analysis.
Final time to reach a steady state, initial and boundary conditions are set as follows:
\begin{alignat*}{4}
	&T = 50, \\
	&h(x,0) = 2 - b(x),\qquad &&h u_m (x,0) = 0  ,\qquad && h \alpha_1 (x,0) = -0.5,\qquad && h \alpha_2 (x,0) = -0.2, \\
	&h(0,t) = 2       ,\qquad &&h u_m (0,t) = 24 ,\qquad && h \alpha_1 (0,t) = -0.5,\qquad && h \alpha_2 (0,t) = -0.2. \\
\end{alignat*}
In particular, we plot $\eta$, $u_m$, $\alpha_1$, and the parabolic velocity profile along the vertical direction computed
at a randomly chosen point $x \approx 18.62$ from:
\begin{equation}\label{eq:veldistribution_parabolic}
	u(\zeta) = u_m + \left(1-2\zeta\right)\alpha_1 + \left(6\zeta^2 - 6\zeta +1 \right)\alpha_2.
\end{equation}

\begin{figure}
	\centering
	\subfigure{
		\includegraphics{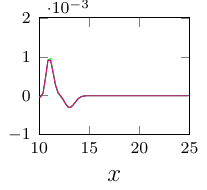}}\quad
	\subfigure{
		\includegraphics{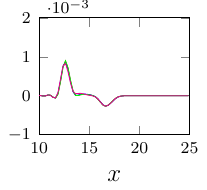}}\quad
	\subfigure{
		\includegraphics{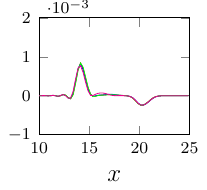}}\quad
	\subfigure{
		\includegraphics{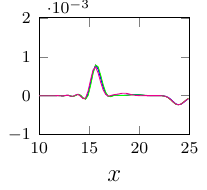}}\quad
	\subfigure{
		\includegraphics{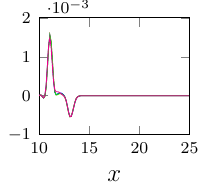}}\quad
	\subfigure{
		\includegraphics{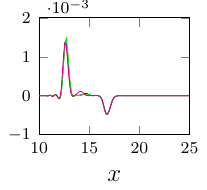}}\quad
	\subfigure{
		\includegraphics{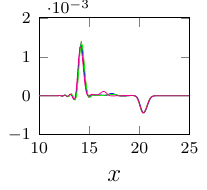}}\quad
	\subfigure{
		\includegraphics{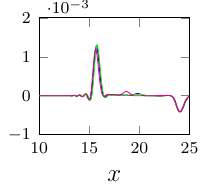}}\quad
	\setcounter{subfigure}{0}
	\subfigure[$t=0.225$]{
		\includegraphics{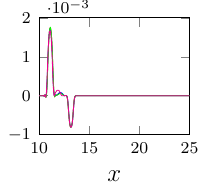}}\quad
	\subfigure[$t=0.45$]{
		\includegraphics{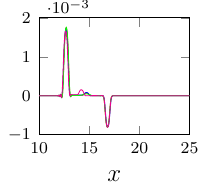}}\quad
	\subfigure[$t=0.675$]{
		\includegraphics{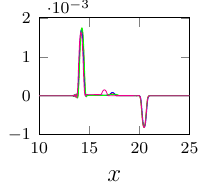}}\quad
	\subfigure[$t=0.9$]{
		\includegraphics{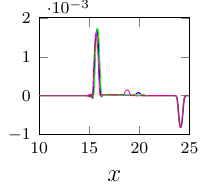}}\quad
	\caption{Small perturbation of supercritical steady state with friction computed with the GF-WENO5 scheme: value of $h-h_{eq}$ for SWME1 (red), SWLME2 (blue), HSWME2 (green) and SWME2 (magenta)  with $N_e=100$ (top), $N_e=200$ (middle), $N_e=800$ (bottom) at different simulation times.}
	\label{fig:super_friction_pert_h}
	\end{figure}

\begin{figure}
		\centering
	\subfigure{
		\includegraphics{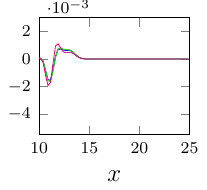}}\quad
	\subfigure{
		\includegraphics{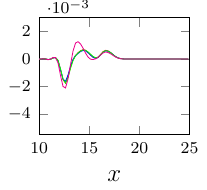}}\quad
	\subfigure{
		\includegraphics{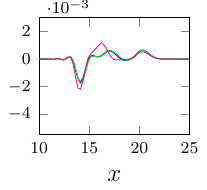}}\quad
	\subfigure{
		\includegraphics{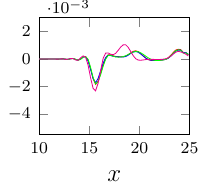}}\quad
	\subfigure{
		\includegraphics{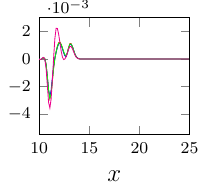}}\quad
	\subfigure{
		\includegraphics{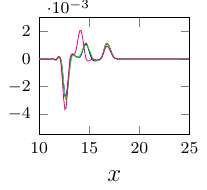}}\quad
	\subfigure{
		\includegraphics{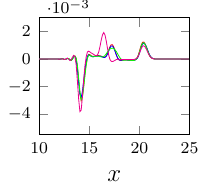}}\quad
	\subfigure{
		\includegraphics{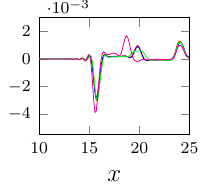}}\quad
	\setcounter{subfigure}{0}
	\subfigure[$t=0.225$]{
		\includegraphics{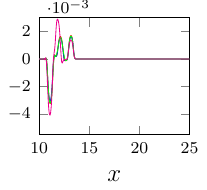}}\quad
	\subfigure[$t=0.45$]{
		\includegraphics{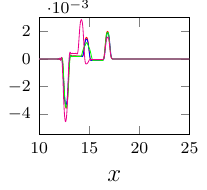}}\quad
	\subfigure[$t=0.675$]{
		\includegraphics{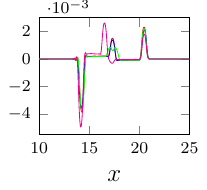}}\quad
	\subfigure[$t=0.9$]{
		\includegraphics{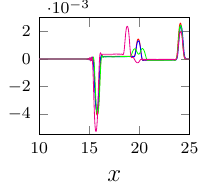}}\quad
	\caption{Small perturbation of supercritical steady state with friction computed with the GF-WENO5 scheme: value of $(h\alpha_1)-(h\alpha_1)_{eq}$ for SWME1 (red), SWLME2 (blue), HSWME2 (green) and SWME2 (magenta)  with $N_e=100$ (top), $N_e=200$ (middle), $N_e=800$ (bottom) at different simulation times.}
	\label{fig:super_friction_pert_ha1}
	\end{figure}

In \cref{fig:super_friction_pert_h} and \cref{fig:super_friction_pert_ha1}, we present the numerical results of the perturbation evolution as deviation from the steady state equilibrium computed with all these models,
using the WENO5 reconstruction on different levels of mesh refinement.
This shows that for a given mesh, the new global flux method can compute the numerical results for the linear models, e.g., SWME1 and SWLME2, and the nonlinear ones HSWME2 and SWME2 with comparable accuracy.
This is promising since for complex nonlinear models like HSWME2 and SWME2, we cannot compute an analytical solution of the steady states, and therefore it allows us to
take as a reference the convergence analysis performed on the simple SWME1 model.

The numerical results of the perturbation analysis allow us to study the differences that are present in the wave propagation of these models.
Firstly, we can notice from the simulations that SWME1 and SWLME2 have the same eigenvalue $u_m$ which results in their superimposed waves.
Interestingly, the HSWME2 presents the two intermediate waves with the double bump, representing the waves with speeds given by $u_m\pm\sqrt{\alpha_1^2/5}$,
while the full SWME2 model provides a different result with waves placed at different positions.
Although the perturbation is rather small, already with a mesh of $N_e=200$ elements the new method is able to capture all the features
in the perturbation analysis, included the double bump waves of the HSWME2 model.

It is also interesting to study the eigenvalues of the system matrices of the different models.
To do so, we take the value of the conservative variables at a specific coordinate position, $x\approx 23$, where all four models
have a similar solution, and compute the eigenvalues of the system matrices with this solution for all considered models.
In \cref{tab:eigenvalues}, we present in descending order both the eigenvalues computed for SWME1, SWLME2 and HSWME2 models,
and the numerical eigenvalues computed with the system matrix of the full SWME2 model.
We show the analytical eigenvalues of the simplified models, as described in \cref{sec:moment}, since they match the numerical eigenvalues of the system matrices,
while the eigenvalues of the full SWME2 model are computed only numerically from the system matrix.
It is interesting to notice that some eigenvalues for SWME2 are very close to some of those presented for simpler models,
while $\lambda_3$ appears to be very different. In particular, this one can be also noticed in the numerical experiments
shown in \cref{fig:super_friction_pert_ha1} where SWME2 presents a wave which is not superimposed in any of the other models.

\begin{table}
	\centering
	\begin{tabular}{c|cccc}\hline\hline
	Eigenvalues  & SWME1 & SWLME2 & HSWME2 & SWME2   \\ \hline
	$\lambda_1$  & 16.15  &  16.19  &  16.15  & 16.17    \\
	$\lambda_2$  & 11.32  &  11.32  &  12.03  & 11.20    \\
	$\lambda_3$  &  --    &   --    &  10.61  & 9.82     \\
	$\lambda_4$  &  6.49  &  6.46   &  6.49   & 6.21     \\ \hline\hline
\end{tabular}
	\caption{Comparison between the analytical eigenvalues of SWME1, SWLME2, HSWME2 models, and the numerical eigenvalues of the full SWME2 model.
	Analytical formulas to compute the eigenvalues of simplified models are taken from \cref{sec:moment}.}\label{tab:eigenvalues}
\end{table}


\section{Conclusion} \label{sec:conclusions}

In this work, we have presented a high order fully well-balanced finite volume method based
on the global flux approach for general hyperbolic systems with non-conservative products,
with the goal of simulating different families shallow water moment models, without a-priori knowledge of analytical steady states.
This is made possible by introducing a tailored quadrature for both source terms and
non-conservative products to include them in a global flux divergence
that allows to rewrite the system in quasi-conservative form.
The steady state preservation property, which is now discretely characterized by a constant
global flux, can then be achieved using suitable numerical fluxes based on the
high order reconstructions of the global flux at interfaces.
The new method is tested using first, third, and fifth order reconstructions to show exact
preservation of lake at rest solutions and the optimal convergence rates on non-trivial moving
equilibria. The impact of reconstruction and mesh refinement is also studied on small
perturbations of the aforementioned steady states proving that the method is
effective for the simulation of complex shallow water moment models.

Future work may include the introduction of additional terms to simulate SWME on manifolds \cite{carlino2023well},
as well as the extension of the method to general two-dimensional models, including two-dimensional shallow water moment models \cite{bauerle2025rotational}.

\section*{Acknowledgments} \label{section: acknowledgments}
This publication is part of the project \textit{HiWAVE} with file number VI.Vidi.233.066 of the \textit{ENW Vidi} research programme, funded by the \textit{Dutch Research Council (NWO)} under the grant \url{https://doi.org/10.61686/CBVAB59929}.

\appendix

\section{Lake at rest exact preservation}\label{app:lar}
To prove the exact preservation of the lake at rest equilibrium we will prove that, for the lake at rest solution, we have
$\bar{\bG}_{i+1}=\bar{\bG}_i$ $\forall\, i$.
Clearly, the first component of the flux concerns the reconstruction of $q\equiv 0$, which is indeed zero.
The same also holds for the velocity moments, i.e.\ $\alpha_1=\ldots=\alpha_N\equiv 0$, involved in the more complex shallow water moment models.
Hence, the only non-zero terms concern the second component of the global flux, meaning the momentum equation of the average velocity.
It is then enough to show that $(\bar\bG_i)_2$ (the second component) is constant to have that $\bar q_i=0$, $\forall\, i$ at the following time steps.
At each quadrature point, we would have
\begin{align}\label{eq:lakeatrest0}
	    (\bar\bG_{i,q})_2 =& (\mathcal R_m)^R_{\iin} + g\frac{(\eta_0-b_{i,q})^2}{2} + g\eta_0\left(\tilde b_{i,q} - b_{\iin}^R\right)  - g\left(\frac{(\tilde b_{i,q})^2}{2} - \frac{(b_{\iin}^R)^2}{2}\right) \\
        =&(\mathcal R_m)^R_{\iin} +g \frac{\eta_0^2}{2} -g\eta_0  b_{\iin}^R +g \frac{(b_{\iin}^R)^2}{2}.
\end{align}
This shows that the global flux is constant across the quadrature points, and thus $(\bar\bG_{i,q})_2 = (\bar\bG_{i})_2$.
We now need to show that this constant is the same $\forall \,i$.
This can be easily shown as follows:
\begin{subequations}
\begin{align*}
	&(\bar\bG_{i+1})_2 - (\bar\bG_{i})_2 \\
	=&  (\mathcal R_m)^R_{\iip} - (\mathcal R_m)^R_{\iin} -g\eta_0  b_{\iip}^R +g \frac{(b_{\iip}^R)^2}{2} +g\eta_0  b_{\iin}^R -g \frac{(b_{\iin}^R)^2}{2}    \\
                =& (\mathcal R_m)^L_{\iip} - (\mathcal R_m)^R_{\iin} + [\![ \mathcal R_m ]\!]_{\iip} -g\eta_0  b_{\iip}^R +g \frac{(b_{\iip}^R)^2}{2} +g\eta_0  b_{\iin}^R -g \frac{(b_{\iin}^R)^2}{2}   \\
                =& \underbrace{g\eta_0\left(b_\iip^L - b_\iin^R\right) - g\left(\frac{(b_\iip^L)^2}{2} - \frac{(b_\iin^R)^2}{2}\right)}_{(\mathcal R_m)^L_{\iip} - (\mathcal R_m)^R_{\iin}} +[\![ \mathcal R_m ]\!]_{\iip}\\
                &\qquad\qquad -g\eta_0  b_{\iip}^R +g \frac{(b_{\iip}^R)^2}{2} +g\eta_0  b_{\iin}^R -g \frac{(b_{\iin}^R)^2}{2} \\
                =& g\eta_0\left(b_\iip^L - b_\iin^R\right) - g\left(\frac{(b_\iip^L)^2}{2} - \frac{(b_\iin^R)^2}{2}\right) +[\![ \mathcal R_m ]\!]_{\iip}=0,\label{eq:last_defi_jump}
\end{align*}
\end{subequations}
recalling that
\begin{equation}
        [\![ \mathcal R_m ]\!]_{\iip} =g\eta_0\left(b_\iip^R - b_\iin^L\right) - g\left(\frac{(b_\iip^R)^2}{2} - \frac{(b_\iin^L)^2}{2}\right),
\end{equation}
which achieves the proof.
\begin{remark}[Definition of the jump of $\mathcal{R}_m$]
	The definition of $[\![ \mathcal R_m ]\!]_{\iip}$ is indeed obtained to achieve well-balancing for lake at rest equilibria.
	In the definition of the jump, the term $\eta_0$ is replaced by its consistent approximation at the interface, i.e.\ $\eta_0 = \frac{1}{2}(\eta_\iip^L + \eta_\iip^R)$.
\end{remark}

 \bibliographystyle{abbrv}
\bibliography{literature}
\end{document}